\documentclass[11pt,a4paper]{article}
\usepackage[margin=27mm]{geometry}
\usepackage{amsmath,amssymb,amsthm,mathtools,microtype}
\usepackage[colorlinks=true,linkcolor=blue,citecolor=blue,urlcolor=blue]{hyperref}
\hypersetup{pdftitle={McKean--Vlasov Diffusions from Liouville Spectral Green Response}}
\allowdisplaybreaks[2]
\numberwithin{equation}{section}
\newtheorem{theorem}{Theorem}[section]
\newtheorem{lemma}[theorem]{Lemma}
\newtheorem{proposition}[theorem]{Proposition}
\newtheorem{corollary}[theorem]{Corollary}
\theoremstyle{remark}
\newtheorem{remark}[theorem]{Remark}
\newcommand{\R}{\mathbb R}
\newcommand{\EE}{\mathbb E}
\newcommand{\PP}{\mathbb P}
\newcommand{\E}{\mathcal E}
\newcommand{\F}{\mathcal F}
\newcommand{\D}{D}
\newcommand{\Law}{\operatorname{Law}}
\newcommand{\one}{\mathbf 1}
\newcommand{\Tr}{\operatorname{Tr}}
\newcommand{\HS}{\mathrm{HS}}
\newcommand{\Ccore}{\mathcal C}
\newcommand{\LL}{\mathcal L}
\title{McKean--Vlasov Diffusions from Liouville Spectral Green Response}
\author{Chunhao Cai\\[3pt]
\small School of Mathematics (Zhuhai), Sun Yat-Sen University, China\\[2pt]
\small \href{mailto:caichh9@mail.sysu.edu.cn}{\texttt{caichh9@mail.sysu.edu.cn}}}
\date{}
\begin{document}
\maketitle
\begin{abstract}
For Dirichlet Liouville Brownian motion with $0<\gamma<2$, the Gaussian
response of a finite vector of logarithmic eigenvalues determines a
conditional Green energy. We weight this energy by the current
spectral density in the fixed reference space. The resulting
self-consistent evolution is a porous-medium equation for the
canonical spectral form. We construct its positive probability flow,
prove an energy identity, and pass from the regularized response forms
to the Green limit. The energy identity yields strong convergence of
the Gaussian gradients and of the conditional response drifts.
A superposition argument then realizes the flow as a McKean--Vlasov
weak diffusion.
\end{abstract}


\section{Introduction and spectral model}
\label{sec:main}
Liouville Brownian motion is the diffusion associated with Liouville
random geometry \cite{GRV}. Its Dirichlet eigenvalues respond to
Cameron--Martin perturbations of the Gaussian field through the
variational calculus of \cite{Cai}. We use this response to define
a Green energy on finite spectral vectors.
In our evolution, the current spectral density weights this energy,
and the corresponding generator evolves that same law. We construct
the resulting nonlinear probability flow and obtain its
McKean--Vlasov weak realization in the Green limit.

Probabilistic representations of irregular porous-medium equations
were studied by Blanchard, R\"ockner and Russo \cite{BRR2010} and,
in the degenerate case, by Barbu, R\"ockner and Russo
\cite{BRR2011}. Barbu and R\"ockner developed this approach for
nonlinear Fokker--Planck equations \cite{BR2018,BR2020} and for
measure initial data \cite{BR2021}. We use the passage from a
nonlinear probability flow to a distribution-dependent SDE described
in \cite[Section 2]{BR2020}. Here the static drift comes from
compensated second spectral response. The Green limit requires
annealed Green-energy bounds and strong convergence of the
density-dependent response, obtained from the nonlinear energy
identity.

\subsection{Spectral variables and Green response}
Let $D\subset\R^2$ be bounded, connected, and smooth, and let $h$ be
the zero-boundary Gaussian free field on $(\Omega,\mathcal A,\PP)$.
We use the Cameron--Martin space
\begin{equation}\label{eq:H}
 H=H_0^1(D),\qquad
 \langle u,v\rangle_H=\frac1{2\pi}\int_D\nabla u\cdot\nabla v\,dx.
\end{equation}
For $0<\gamma<2$, write $M_h=M_\gamma(h)$ for the associated Gaussian
multiplicative chaos measure with full-variance normalization. Let
$0<\lambda_1<\lambda_2<\cdots$ be the almost surely simple Dirichlet
Liouville eigenvalues, with real eigenfunctions normalized by
\begin{equation}\label{eq:spectrum}
 \langle\phi_i,v\rangle_H
   =\lambda_i\int_D\phi_i v\,dM_h,\qquad
 \int_D\phi_i\phi_j\,dM_h=\delta_{ij},\qquad v\in H.
\end{equation}
Thus $\nu_i=\phi_i^2M_h$ is a probability measure. Put
$Y_i=\log\lambda_i$ for $i\ge1$, fix $N<\infty$, and set
$Y=(Y_1,\ldots,Y_N)$. Throughout, $\rho=\Law(Y)$ is the reference
measure and $E=\operatorname{supp}\rho\subset\R^N$ is the state space.

Let $G_D(x,z)$ denote the Dirichlet Green kernel corresponding to
the normalization in \eqref{eq:H}.
For a measure $\eta$ of finite Green energy, its Green--Riesz
potential $U_D\eta\in H$ is characterized by
\[
 \langle U_D\eta,g\rangle_H=\int_Dg\,d\eta,\qquad g\in H.
\]
Write $I_G(\eta,\zeta)=\langle U_D\eta,U_D\zeta\rangle_H$ and
$I_G(\eta)=I_G(\eta,\eta)$. With $\D$ denoting the Malliavin
derivative of the Gaussian field, the spectral response is
\begin{equation}\label{eq:response}
 \D Y_i=-\gamma U_D\nu_i.
\end{equation}
The Green-energy estimates in Appendix~\ref{app:annealed} justify
this identity in $\mathbb D^{1,r}$ for every $1<r<\infty$;
the gradient identification is proved in Appendix~\ref{app:forms}.
Conditioning the scalar products of these responses on $Y$ gives
the matrix $a=(a_{ij})_{i,j\le N}$, with
\begin{equation}\label{eq:green-covariance}
 a_{ij}(y)=\gamma^2\EE[I_G(\nu_i,\nu_j)\mid Y=y].
\end{equation}
We write $\sigma=(2a)^{1/2}$ for its symmetric nonnegative square
root.

\subsection{Response energies and the nonlinear evolution}
Let $-\Delta_D$ be the Dirichlet Laplacian and set
$Q_s=(I-\Delta_D)^{-s}$ for $s\ge0$. On smooth spectral functions,
the regularized response energy is
\begin{equation}\label{eq:intro-response-energy}
 \E_s(u,v)=\EE\langle Q_s\D(u(Y)),\D(v(Y))\rangle_H,
 \qquad u,v\in C_c^\infty(\R^N)|_E.
\end{equation}
For $s>0$, let $(\E_s,\F_s)$ be its closure, with closed response
gradient $D_s$. We use the canonical increasing limit
$(\E_*,\F_*)$ as $s\downarrow0$. The domains and the construction
of $D_s$ are given in Section~\ref{sec:law}. On $\F_*$, the
limiting energy has the same Gaussian representation with $Q_0=I$.

The static drift $b=(b_i)\in L^1(\rho;\R^N)$ is the weak divergence
of $a$ with respect to $\rho$. Appendix~\ref{app:drift} constructs it
from the compensated second spectral response. The conditional
responses on the two form domains are
\[
 \begin{aligned}
 R_{i,s}[g](y)
 &=\EE[\langle Q_s^{1/2}\D Y_i,D_sg\rangle_H\mid Y=y],
       &&g\in\F_s,\\
 R_i[g](y)
 &=\EE[\langle\D Y_i,\D(g(Y))\rangle_H\mid Y=y],
       &&g\in\F_*.
 \end{aligned}
\]
Write $R_s[g]=(R_{i,s}[g])_{i\le N}$ and $R[g]=(R_i[g])_{i\le N}$.
All conditional coefficients are understood $\rho$-almost everywhere.

For a current law $\mu=f\rho$, we weight the response energy by $f$
in the fixed space $L^2(\rho)$. On smooth functions this energy is
$\int_E f\,\nabla u^Ta\nabla v\,d\rho$.
For a bounded density $f\in\F_*$, the product rule gives diffusion
matrix $fa$ and drift $fb+R[f]$; Section~\ref{sec:law} proves this
generator identity. Self-consistency pairs the weighted energy with
the current density itself. The chain rule
$f(Y)\D(f(Y))=\frac12\D((f^2)(Y))$ then yields the quadratic
evolution stated below.

\subsection{Main results}
Fix $0<c\le1\le C$ and an initial density satisfying
\begin{equation}\label{eq:init}
 f_0\in L^\infty(\rho),\qquad c\le f_0\le C,\qquad
 \int_E f_0\,d\rho=1.
\end{equation}
The parameter $s>0$ regularizes the Gaussian response, while $t$
denotes time. We use $s=*$ for the canonical Green form and write
$f_{*,t}=f_t$. The density equation is
\begin{equation}\label{eq:nlweak}
 \int_Eu f_{s,t}\,d\rho-\int_Eu f_0\,d\rho
 =-\frac12\int_0^t\E_s(f_{s,r}^2,u)\,dr,
 \qquad u\in\F_s.
\end{equation}

\begin{theorem}\label{thm:flow}
For each $s>0$ and for $s=*$, equation \eqref{eq:nlweak} has a
solution $f_s$ such that, for every $T<\infty$,
$f_s\in C([0,T];L^1(\rho))$, $c\le f_{s,t}\le C$, and
$\int_E f_{s,t}\,d\rho=1$ for $0\le t\le T$.
Moreover, $f_s^2\in L^2(0,T;\F_s)$ and
$\partial_t f_s\in L^2(0,T;\F_s')$, and
\begin{equation}\label{eq:thmflow3}
 \int_E\frac{f_{s,T}^3}{3}\,d\rho
 +\frac12\int_0^T\E_s(f_{s,t}^2,f_{s,t}^2)\,dt
 =\int_E\frac{f_0^3}{3}\,d\rho.
\end{equation}
\end{theorem}

Density convergence and \eqref{eq:thmflow3} give convergence of the
energy norms of $f_s^2$. Combined with weak gradient identification,
this yields strong convergence of $D_s(f_{s,t}^2)$. The bound
$f_{s,t}\ge c$ then transfers this convergence to $D_sf_{s,t}$ and
to the conditional response entering the drift.

\begin{theorem}\label{thm:limit}
Let $f_s$ be the flows constructed in the proof of
Theorem~\ref{thm:flow}. For every
$T<\infty$, as $s\downarrow0$,
\begin{align}
 \sup_{0\le t\le T}\|f_{s,t}-f_t\|_{L^1(\rho)}
 &\longrightarrow0,\label{eq:thmlimit1}\\
 \int_0^T\|R_{i,s}[f_{s,t}]-R_i[f_t]\|_{L^1(\rho)}\,dt
 &\longrightarrow0,\qquad 1\le i\le N.\label{eq:thmlimit2}
\end{align}
The limiting flow is realized on a filtered probability space by
a continuous $E$-valued process $X$ and an $N$-dimensional Brownian
motion $W$ such that $\Law(X_t)=f_t\rho$ for $0\le t\le T$ and
\begin{equation}\label{eq:thmsde}
 \begin{aligned}
 dX_i(t)
 &=\bigl[f_t(X_t)b_i(X_t)+R_i[f_t](X_t)\bigr]\,dt\\
 &\quad+\sqrt{f_t(X_t)}
               \sum_{j=1}^N\sigma_{ij}(X_t)\,dW_j(t),
 \qquad 1\le i\le N.
 \end{aligned}
\end{equation}
\end{theorem}

Section~\ref{sec:law} derives the nonlinear equation from the weighted
response energy. Section~\ref{sec:flow} constructs its probability
flow and proves Theorem~\ref{thm:flow}. The estimates of
Section~\ref{sec:green} are used in the final proof of
Theorem~\ref{thm:limit} in Section~\ref{sec:realize}.
The static spectral analysis is given in Appendix~\ref{sec:appendix};
Appendix~\ref{app:time-chain} proves the time-chain rule used in
the nonlinear energy identity.

\section{From Green response to law-dependent evolution}
\label{sec:law}
The evolution is built in the fixed space $L^2(\rho)$. We first
specify the closed response forms and their integration-by-parts
identities, then use the current density as their weight.

\subsection{Closed response forms}
Let $-\Delta_D\psi_k=\kappa_k\psi_k$, with $(\psi_k)$ orthonormal
in $L^2(D)$, and write
\begin{equation}\label{eq:Q}
 e_k=\sqrt{\frac{2\pi}{\kappa_k}}\,\psi_k,\qquad
 Q_se_k=q_k(s)e_k,\qquad q_k(s)=(1+\kappa_k)^{-s},\quad s\ge0.
\end{equation}
The $(e_k)$ form an orthonormal basis of $H$, and $Q_0=I$.
The Gaussian response matrix $A_s=(A_{ij,s})_{i,j\le N}$ has
entries $A_{ij,s}=\langle Q_s\D Y_i,\D Y_j\rangle_H$.
Its conditional version is
\begin{equation}\label{eq:a}
 a_{ij,s}(y)=\EE[A_{ij,s}\mid Y=y],
 \qquad a_s=(a_{ij,s})_{i,j\le N}.
\end{equation}
Thus $a_0=a$ by \eqref{eq:response}. On
$\Ccore=C_c^\infty(\R^N)|_E$, the chain rule
$\D(u(Y))=\sum_i(\partial_i u)(Y)\D Y_i$ converts
\eqref{eq:intro-response-energy} into
\begin{equation}\label{eq:form_s}
 \begin{aligned}
 \E_s(u,v)
 &=\sum_{i,j=1}^N\EE\bigl[
       (\partial_i u)(Y)(\partial_j v)(Y)
       A_{ij,s}\bigr]\\
 &=\int_E\nabla u^Ta_s\nabla v\,d\rho.
 \end{aligned}
\end{equation}

Let $D_s$ be the closure, from $L^2(\rho)$ to $L^2(\PP;H)$, of
$u\mapsto Q_s^{1/2}\D(u(Y))$ on $\Ccore$. For $s>0$, its domain
is $\F_s$, the minimal form domain of $\E_s$; at $s=0$ we write
$(\E_0^{\min},\F_0^{\min})$ for the minimal closure.
As $s\downarrow0$, the weights in \eqref{eq:Q} increase, and the
canonical increasing form is
\begin{equation}\label{eq:canonical-intro}
 \F_*=\left\{u\in\bigcap_{s>0}\F_s:
                 \sup_{s>0}\E_s(u,u)<\infty\right\},\qquad
 \E_*(u,u)=\sup_{s>0}\E_s(u,u).
\end{equation}
The compensated coefficients $b_s=(b_{i,s})$ are constructed in
Appendix~\ref{app:drift}, with $b_0=b$. They enter the differential
expressions
\begin{equation}\label{eq:static_generator}
 \ell_su=a_s:\nabla^2u+b_s\cdot\nabla u,\qquad
 \ell u=a:\nabla^2u+b\cdot\nabla u,\qquad u\in C_c^\infty(\R^N).
\end{equation}
Here $\nabla^2u=(\partial_{ij}u)_{i,j\le N}$ is the Euclidean Hessian.
The form identities below specify how these expressions act under
integration against $\rho$.

\begin{proposition}\label{prop:static}
For $s>0$, $(\E_s,\F_s)$ is a regular, strongly local,
conservative Dirichlet form, with $1,y_i\in\F_s$ and
\begin{equation}\label{eq:static-gradients}
 \E_s(u,v)=\EE\langle D_su,D_sv\rangle_H,\qquad
 D_sy_i=Q_s^{1/2}\D Y_i.
\end{equation}
The forms $\E_s$ Mosco-converge to $\E_*$ as $s\downarrow0$.
The canonical form is closed, strongly local, and conservative.
Moreover, $\F_0^{\min}\subseteq\F_*$, $1,y_i\in\F_*$, and
$u(Y)\in\mathbb D^{1,2}$ for $u\in\F_*$, with
\begin{equation}\label{eq:static-canonical-gradient}
 \E_*(u,v)=\EE\langle\D(u(Y)),\D(v(Y))\rangle_H,
 \qquad u,v\in\F_*.
\end{equation}
The coefficients satisfy $a_{ij,s},a_{ij}\in L^r(\rho)$ for
every finite $r\ge1$ and $b_{i,s},b_i\in L^1(\rho)$, and
\begin{equation}\label{eq:static-coefficients}
 \|a_s-a\|_{L^r(\rho)}\longrightarrow0\quad(1\le r<\infty),
 \qquad \|b_s-b\|_{L^1(\rho)}\longrightarrow0
 \qquad(s\downarrow0).
\end{equation}
For $u\in C_c^\infty(\R^N)$, the weak generator identities are
\begin{equation}\label{eq:static-form-identities}
 \begin{aligned}
 \E_s(u,g)&=-\int_E g\,\ell_su\,d\rho,
       &&g\in\F_s\cap L^\infty(\rho),\\
 \E_*(u,g)&=-\int_E g\,\ell u\,d\rho,
       &&g\in\F_*\cap L^\infty(\rho).
 \end{aligned}
\end{equation}
\end{proposition}
\begin{proof}
Proposition~\ref{prop:forms} proves closability and gives the
regularized forms, their gradients, and the coordinate functions.
Lemma~\ref{lem:weighted-gradient-domain} and
Theorem~\ref{thm:Alimit} identify the increasing limit and its
Gaussian gradient; the inclusion of the minimal domain also gives
$1,y_i\in\F_*$. The coefficient estimates follow from
Lemma~\ref{lem:MalliavinY} and Proposition~\ref{prop:Acoeff}.
Proposition~\ref{prop:IBP} and
Lemma~\ref{lem:appendix-endpoint-ibp} give
\eqref{eq:static-form-identities}.
\end{proof}

For the weighted construction, it is convenient to write the energies
as integrals of their conditional densities:
\begin{align}
 \Gamma_s(u,v)(y)
 &=\EE[\langle D_su,D_sv\rangle_H\mid Y=y],
       &&u,v\in\F_s,\label{eq:gamma-s}\\
 \Gamma(u,v)(y)
 &=\EE[\langle\D(u(Y)),\D(v(Y))\rangle_H\mid Y=y],
       &&u,v\in\F_*.\label{eq:gamma}
\end{align}
The tower property gives $\E_s(u,v)=\int_E\Gamma_s(u,v)\,d\rho$
and $\E_*(u,v)=\int_E\Gamma(u,v)\,d\rho$.
Cauchy--Schwarz gives
$\|\Gamma(u,v)\|_{L^1(\rho)}
 \le\E_*(u,u)^{1/2}\E_*(v,v)^{1/2}$, and the same bound holds for
$\Gamma_s$ and $\E_s$.
By \eqref{eq:static-gradients} and
\eqref{eq:static-canonical-gradient}, the responses of Section~\ref{sec:main}
can therefore be written as
\begin{equation}\label{eq:R}
 R_{i,s}[g]=\Gamma_s(y_i,g),\quad g\in\F_s,\qquad
 R_i[g]=\Gamma(y_i,g),\quad g\in\F_*.
\end{equation}

\subsection{The law-dependent generator}
Let $\mu=f\rho$ be a probability law with $c\le f\le C$.
The original spectral law $\rho=\Law(Y)$ and its response energy
$\Gamma$ remain fixed; law dependence enters by weighting this energy
with the current density $f=d\mu/d\rho$ in $L^2(\rho)$:
\begin{equation}\label{eq:weighted}
 \E_\mu(u,v)=\int_E f\,\Gamma(u,v)\,d\rho,\qquad u,v\in\F_*.
\end{equation}
The bounds on $f$ give
\begin{equation}\label{eq:weightedclose}
 c\,\E_*(u,u)\le\E_\mu(u,u)\le C\,\E_*(u,u),
 \qquad u\in\F_*,
\end{equation}
so the form norm is equivalent to that of $\E_*$. Hence
$\E_\mu$ is closed. The chain rule makes it Markovian and strongly
local, while $1\in\F_*$ and $\E_\mu(1,1)=0$ give conservativeness.
To identify its drift, we use the product rule with the density $f$.

\begin{lemma}\label{lem:lawgen}
If $f\in\F_*$ and $c\le f\le C$, then for
$u\in C_c^\infty(\R^N)$ the expression
$\ell[f]u=f\,\ell u+\Gamma(f,u)$ belongs to $L^1(\rho)$ and satisfies
\begin{equation}\label{eq:weightedgen}
 -\E_\mu(u,g)=\int_Eg\,\ell[f]u\,d\rho,\qquad
 g\in\F_*\cap L^\infty(\rho).
\end{equation}
Its coefficients are
\begin{equation}\label{eq:law-coefficients}
 a[f]=f a,\qquad b_i[f]=f b_i+R_i[f],\qquad
 \ell[f]u=a[f]:\nabla^2u+b[f]\cdot\nabla u.
\end{equation}
\end{lemma}
\begin{proof}
For bounded $f,g\in\F_*$, the regularized product estimate
\cite{FOT} gives
\[
 \E_s(fg,fg)\le
 2\|f\|_\infty^2\E_s(g,g)+2\|g\|_\infty^2\E_s(f,f).
\]
Thus $fg\in\F_s$ for every $s>0$, and taking the supremum over $s$
gives $fg\in\F_*$. We may use this product in
\eqref{eq:static-form-identities}. The Gaussian product rule yields
\begin{equation}\label{eq:product}
 \begin{aligned}
 -\E_\mu(u,g)
 &=-\E_*(u,fg)+\int_Eg\,\Gamma(u,f)\,d\rho\\
 &=\int_Eg\bigl(f\,\ell u+\Gamma(f,u)\bigr)\,d\rho.
 \end{aligned}
\end{equation}
The coefficient bounds of Proposition~\ref{prop:static} and the
$L^1$ bound on $\Gamma$ show that the integrand defining
$\ell[f]u$ is in $L^1(\rho)$. This proves \eqref{eq:weightedgen}.
Finally, applying the chain rule to $u(Y)$ gives
\begin{equation}\label{eq:response-chain}
 \Gamma(f,u)=\sum_{i=1}^N(\partial_i u)R_i[f].
\end{equation}
Substitution of \eqref{eq:static_generator} identifies the coefficients
in \eqref{eq:law-coefficients}.
\end{proof}

For the regularized form
$\E_{s,\mu}(u,v)=\int_E f\,\Gamma_s(u,v)\,d\rho$, the same
calculation applies when $f\in\F_s$. It gives
$\ell_s[f]u=f\,\ell_su+\Gamma_s(f,u)$ and
\begin{equation}\label{eq:regularized-law-coefficients}
 a_s[f]=f a_s,\qquad b_s[f]=f b_s+R_s[f].
\end{equation}

\subsection{The self-consistent equation}
Let $L_*$ be the nonpositive $L^2(\rho)$ generator of the canonical
form. We use the same symbol for its extension to the form dual,
characterized by
\begin{equation}\label{eq:form-generator-duality}
 \langle L_*v,u\rangle_{\F_*',\F_*}=-\E_*(v,u),
 \qquad u,v\in\F_*.
\end{equation}
Here $\F_*'$ is the dual of $\F_*$ with pivot space $L^2(\rho)$;
the notation $\F_s'$ has the analogous meaning.

Let $f\in L^2(0,T;\F_*)$ be a density curve with $c\le f_t\le C$,
and put $\mu_t=f_t\rho$. We formulate self-consistency on the full
canonical domain:
\begin{equation}\label{eq:forward-form}
 \frac{d}{dt}\int_Eu f_t\,d\rho=-\E_{\mu_t}(u,f_t),
 \qquad u\in\F_*,
\end{equation}
in the sense of distributions in time. For smooth $u$, the right-hand
side is $\int_E\ell[f_t]u\,d\mu_t$ by Lemma~\ref{lem:lawgen}.
The bounded chain rule gives $f^2\in L^2(0,T;\F_*)$ and
\begin{equation}\label{eq:self-consistency}
 \begin{aligned}
 \E_{\mu_t}(u,f_t)
 &=\EE\bigl[f_t(Y)
       \langle\D(u(Y)),\D(f_t(Y))\rangle_H\bigr]\\
 &=\frac12\EE\langle\D(u(Y)),\D((f_t^2)(Y))\rangle_H\\
 &=\frac12\E_*(u,f_t^2).
 \end{aligned}
\end{equation}
Combining \eqref{eq:forward-form} with
\eqref{eq:form-generator-duality} therefore gives
\begin{equation}\label{eq:porous-medium}
 \partial_t f_t=\frac12L_*(f_t^2)
 \quad\text{in }L^2(0,T;\F_*').
\end{equation}
With initial density $f_0$, its time-integrated form is
\eqref{eq:nlweak} at $s=*$.

For $s>0$, apply the same construction to a density curve
$f_s\in L^2(0,T;\F_s)$ with the same bounds. Its chain rule reads
$D_s(f_{s,t}^2)=2f_{s,t}(Y)D_sf_{s,t}$, and the resulting
time-integrated equation is \eqref{eq:nlweak}.
Conservativeness makes the constant density $f_t=1$ an equilibrium.
We now construct the nonlinear flow for the initial data
\eqref{eq:init} and prove its cubic energy identity.

\section{Nonlinear probability flow}\label{sec:flow}
Fix one of the forms $(\E_s,\F_s)$, $s>0$, or $(\E_*,\F_*)$.
In this section write $(\E,\F)$ for this form and $L,P_t$ for its
generator and linear semigroup, with the duality convention of
\eqref{eq:form-generator-duality}.
The constants $c,C$ and the initial density are as in \eqref{eq:init}.
We work on the closed order interval
$K=\{u\in L^1(\rho):c\le u\le C\}$ and equip $\F$ with the
form norm $\|u\|_{\F}^2=\|u\|_2^2+\E(u,u)$.

\subsection{The nonlinear resolvent}
Writing $v=u^2$ makes the elliptic term in one implicit time step
linear. Choose a continuously differentiable increasing function
$\psi:\R\to\R$ equal to $\sqrt r$ on $[c^2,C^2]$, with
$0<m\le\psi'\le M$, and set $\Psi(r)=\int_0^r\psi(q)\,dq$.

\begin{lemma}\label{lem:resolvent}
For $\tau>0$ and $g\in K$, there is a unique $u\in K$ with
$u^2\in\F$ satisfying
\begin{equation}\label{eq:res}
 \int_E(u-g)w\,d\rho+\frac\tau2\E(u^2,w)=0,\qquad w\in\F.
\end{equation}
The solution map $J_\tau:g\mapsto u$ preserves the integral and
satisfies, for $g_1,g_2\in K$,
\begin{equation}\label{eq:rescontr}
 \begin{aligned}
 \int_E(J_\tau g_1-J_\tau g_2)^+\,d\rho
 &\le\int_E(g_1-g_2)^+\,d\rho,\\
 \|J_\tau g_1-J_\tau g_2\|_1&\le\|g_1-g_2\|_1.
 \end{aligned}
\end{equation}
\end{lemma}

\begin{proof}
On $\F$, consider the functional
\begin{equation}\label{eq:J}
 F_{\tau,g}(v)=\int_E\bigl(\Psi(v)-gv\bigr)d\rho+\frac\tau4\E(v,v).
\end{equation}
The bound $\psi'\ge m$ gives strict convexity and
\begin{equation}\label{eq:Jcoerc}
 F_{\tau,g}(v)\ge\frac m4\|v\|_2^2
       +\frac\tau4\E(v,v)-\frac1m\|g-\psi(0)\|_2^2.
\end{equation}
Coercivity in the form norm and weak lower semicontinuity therefore
give a unique minimizer $v\in\F$. Its Euler equation is
\begin{equation}\label{eq:euler}
 \int_E(\psi(v)-g)w\,d\rho+\frac\tau2\E(v,w)=0,\qquad w\in\F.
\end{equation}
The truncation $(v-C^2)^+$ belongs to $\F$. Testing with it and
using $g\le C$ yields
\[
 0\le\int_E(\psi(v)-C)(v-C^2)^+\,d\rho
 \le-\frac\tau2\E(v,(v-C^2)^+)\le0.
\]
Strict monotonicity of $\psi$ gives $v\le C^2$. Testing with
$-(c^2-v)^+$ gives $v\ge c^2$ in the same way. Thus
$u=\psi(v)=\sqrt v\in K$, and \eqref{eq:euler} becomes
\eqref{eq:res}.

To compare two solutions, put $v_j=u_j^2$ and
$\eta_\varepsilon(r)=\min\{1,r^+/\varepsilon\}$. The Dirichlet
form contraction rule gives
$\E(v_1-v_2,\eta_\varepsilon(v_1-v_2))\ge0$.
Subtracting their equations and testing with
$\eta_\varepsilon(v_1-v_2)$ therefore gives
\[
 \begin{aligned}
 \int_E(u_1-u_2)\eta_\varepsilon(v_1-v_2)\,d\rho
 &\le\int_E(g_1-g_2)\eta_\varepsilon(v_1-v_2)\,d\rho\\
 &\le\int_E(g_1-g_2)^+\,d\rho.
 \end{aligned}
\]
Since $u_1-u_2$ and $v_1-v_2$ have the same sign, dominated
convergence as $\varepsilon\downarrow0$ proves the first inequality
in \eqref{eq:rescontr}. Exchanging the indices and adding proves
the second inequality and uniqueness. Finally, $1\in\F$ and
$\E(v,1)=0$, so the test $w=1$ in \eqref{eq:res} proves preservation
of the integral.
\end{proof}

We next identify the operator whose resolvents have just been
constructed. Set
\begin{equation}\label{eq:nonlinear-op}
 \operatorname{Dom}(A)=\{u\in K:u^2\in \operatorname{Dom}(L)\},\qquad
 Au=-\tfrac12 L(u^2).
\end{equation}
Here $\operatorname{Dom}(L)$ is the $L^2(\rho)$ generator domain. For
$u=J_\tau g$, the boundedness of $u-g$ and \eqref{eq:res}
give $u^2\in \operatorname{Dom}(L)$, with
\begin{equation}\label{eq:graph-res}
 L(u^2)=\frac2\tau(u-g),\qquad (I+\tau A)J_\tau g=g.
\end{equation}
For $u_1,u_2\in \operatorname{Dom}(A)$, write
$g_j=u_j+\tau Au_j\in L^2(\rho)$. The same positive-part test
applies to their weak equations. Exchanging the indices and adding
now yields
\begin{equation}\label{eq:accretivity}
 \|u_1-u_2\|_{L^1(\rho)}
 \le\|u_1-u_2+\tau(Au_1-Au_2)\|_{L^1(\rho)},\qquad \tau>0.
\end{equation}
Thus $A$ is accretive in $L^1(\rho)$. For $0<\lambda<\tau$,
rewriting $(I+\tau A)u=g$ with step $\lambda$ gives the datum
$(\lambda/\tau)g+(1-\lambda/\tau)u\in K$. Uniqueness in
Lemma~\ref{lem:resolvent} then gives the resolvent identity
\begin{equation}\label{eq:residentity}
 J_\tau g=J_\lambda\left(\frac\lambda\tau g+
 \left(1-\frac\lambda\tau\right)J_\tau g\right).
\end{equation}

It remains to check density of $\operatorname{Dom}(A)$ in $K$. First let $g\in K$
with $g^2\in\F$, and put $v_\tau=(J_\tau g)^2$.
Since $\Psi'(g^2)=g$, strong convexity of $\Psi$ and comparison
of \eqref{eq:J} at $v_\tau$ and $g^2$ yield
\begin{equation}\label{eq:resid}
 \begin{aligned}
 \frac m2\|v_\tau-g^2\|_2^2+\frac\tau4\E(v_\tau,v_\tau)
 &\le\frac\tau4\E(g^2,g^2),\\
 \|J_\tau g-g\|_2&\le\frac1{2c}\|v_\tau-g^2\|_2
       \longrightarrow0\qquad(\tau\downarrow0).
 \end{aligned}
\end{equation}
For general $g\in K$, use $g_\varepsilon=P_\varepsilon g$.
The conservative Markov semigroup gives $g_\varepsilon\in K\cap\F$;
the bounded chain rule gives $g_\varepsilon^2\in\F$.
Contractivity of $J_\tau$ and \eqref{eq:resid} imply
\begin{equation}\label{eq:densityres}
 \begin{aligned}
 \|J_\tau g-g\|_1
 &\le2\|g-g_\varepsilon\|_1+\|J_\tau g_\varepsilon-g_\varepsilon\|_1,\\
 \limsup_{\tau\downarrow0}\|J_\tau g-g\|_1
 &\le2\|g-P_\varepsilon g\|_1\longrightarrow0
       \qquad(\varepsilon\downarrow0).
 \end{aligned}
\end{equation}
Since $J_\tau g\in \operatorname{Dom}(A)$ by \eqref{eq:graph-res}, this proves
the domain density. Together with the resolvent construction, we have
\begin{equation}\label{eq:nonlinear-range}
 \overline{\operatorname{Dom}(A)}^{\,L^1(\rho)}=K,\qquad
 K\subset\operatorname{Ran}(I+\tau A)\quad(\tau>0).
\end{equation}

\subsection{Final proof of Theorem~\ref{thm:flow}}
\label{subsec:final-flow}
Fix $0<T<\infty$. For $\tau>0$, start the implicit Euler scheme
from the prescribed initial density:
\begin{equation}\label{eq:euler-scheme}
 f_\tau^0=f_0,\qquad f_\tau^k=J_\tau f_\tau^{k-1},
 \qquad k\ge1.
\end{equation}
Lemma~\ref{lem:resolvent} gives
\begin{equation}\label{eq:discrete-probability}
 c\le f_\tau^k\le C,\qquad
 \int_E f_\tau^k\,d\rho=1,
\end{equation}
and each step satisfies
\begin{equation}\label{eq:euler-step}
 \int_E(f_\tau^k-f_\tau^{k-1})w\,d\rho
 +\frac{\tau}{2}\E((f_\tau^k)^2,w)=0,\qquad w\in\F.
\end{equation}
Convexity of $r\mapsto r^3/3$ on $[c,C]$ gives
$(a-b)a^2\ge(a^3-b^3)/3$. Testing \eqref{eq:euler-step}
with $(f_\tau^k)^2$ and summing from $k=1$ to $n$ gives
\begin{equation}\label{eq:discrete-ent}
 \int_E\frac{(f_\tau^n)^3}{3}\,d\rho
 +\frac{\tau}{2}\sum_{k=1}^n
       \E((f_\tau^k)^2,(f_\tau^k)^2)
 \le\int_E\frac{f_0^3}{3}\,d\rho.
\end{equation}

The operator in \eqref{eq:nonlinear-op} is accretive by
\eqref{eq:accretivity}, its domain is dense in $K$, and its
resolvent range contains $K$ by \eqref{eq:nonlinear-range}.
The Crandall--Liggett theorem \cite{CL} therefore yields
a contraction semigroup on $K$. Its orbit from $f_0$ satisfies
$f\in C([0,T];L^1(\rho))$ and, for the piecewise constant
interpolation of the Euler values,
\begin{equation}\label{eq:eulerconv}
 \sup_{0\le t\le T}
   \|\bar f_\tau(t)-f_t\|_{L^1(\rho)}\longrightarrow0,
 \qquad
 \bar f_\tau(t)=f_\tau^k\quad((k-1)\tau<t\le k\tau),
\end{equation}
with $\bar f_\tau(0)=f_0$. Passing to the limit in
\eqref{eq:discrete-probability} gives $c\le f_t\le C$ and
$\int_E f_t\,d\rho=1$ for every $t$.

To pass the discrete equation to the limit, let
$\widetilde f_\tau$ be the affine interpolation through
$\widetilde f_\tau(k\tau)=f_\tau^k$, and put
$\bar v_\tau=\bar f_\tau^2$. Repeated use of
\eqref{eq:rescontr} bounds every Euler increment by the first one,
so
\[
 \sup_{0\le t\le T}
       \|\widetilde f_\tau(t)-\bar f_\tau(t)\|_1
 \le\|J_\tau f_0-f_0\|_1\longrightarrow0
\]
by \eqref{eq:densityres}. Both interpolations thus converge to
$f$. Equation~\eqref{eq:euler-step} also gives
\begin{equation}\label{eq:euler-dual}
 \partial_t\widetilde f_\tau=\frac12L\bar v_\tau
 \quad\text{in }\F',\qquad
 \int_0^T\|\partial_t\widetilde f_\tau\|_{\F'}^2\,dt
 \le\frac14\int_0^T\E(\bar v_\tau,\bar v_\tau)\,dt.
\end{equation}
The discrete energy bound \eqref{eq:discrete-ent} makes
$\bar v_\tau$ bounded in $L^2(0,T;\F)$ and
$\partial_t\widetilde f_\tau$ bounded in $L^2(0,T;\F')$.
For $a,b\in[c,C]$, we have $|a-b|^2\le2C|a-b|$ and
$|a^2-b^2|^2\le8C^3|a-b|$. Thus the $L^1$ convergence and
the weak energy bound identify the limits as
\begin{equation}\label{eq:squareweak}
 \begin{aligned}
 \bar v_\tau&\rightharpoonup f^2
       &&\text{in }L^2(0,T;\F),\\
 \bar v_\tau&\longrightarrow f^2
       &&\text{in }L^2(0,T;L^2(\rho)),\\
 \widetilde f_\tau&\longrightarrow f
       &&\text{in }L^2(0,T;L^2(\rho)).
 \end{aligned}
\end{equation}
For $w\in\F$ and $\chi\in C_c^\infty(0,T)$, testing
\eqref{eq:euler-dual} and using \eqref{eq:squareweak} now yields
\begin{equation}\label{eq:timeweak}
 \begin{aligned}
 -\int_0^T\chi'(t)\int_E f_t w\,d\rho\,dt
 &=\lim_{\tau\downarrow0}\int_0^T\chi(t)
       \langle\partial_t\widetilde f_\tau,w\rangle_{\F',\F}\,dt\\
 &=-\frac12\int_0^T\chi(t)\E(f_t^2,w)\,dt.
 \end{aligned}
\end{equation}
Hence $\partial_t f=\tfrac12L(f^2)\in L^2(0,T;\F')$.
The density bound and $L^1$ continuity imply
$\|f_t-f_r\|_2^2\le2C\|f_t-f_r\|_1$, so
$f\in C([0,T];L^2(\rho))$. Integrating \eqref{eq:timeweak}
with the initial trace $f_0$ proves \eqref{eq:nlweak} for every
$w\in\F$ and every $t\in[0,T]$.

It remains to obtain the energy identity. The chain rule for
$r\mapsto\sqrt r$ on $[c^2,C^2]$ gives
\[
 \int_0^T\E(f_t,f_t)\,dt
 \le\frac1{4c^2}\int_0^T\E(f_t^2,f_t^2)\,dt<\infty.
\]
Thus $f\in L^2(0,T;\F)$, and Lemma~\ref{lem:time-chain} applies.
Using the limiting equation in its dual pairing gives
\begin{equation}\label{eq:energy-equal}
 \begin{aligned}
 \int_E\frac{f_T^3-f_0^3}{3}\,d\rho
 &=\int_0^T\langle\partial_t f_t,f_t^2\rangle_{\F',\F}\,dt\\
 &=-\frac12\int_0^T\E(f_t^2,f_t^2)\,dt.
 \end{aligned}
\end{equation}
This proves \eqref{eq:thmflow3}. The construction applies to every
regularized form and to the canonical form, giving $f_{s,t}$ and
$f_t$ with the asserted bounds and regularity.

\section{Estimates for the nonlinear Green limit}
\label{sec:green}
We prepare the resolvent, Gaussian gradient, and conditional
response estimates used to pass from the regularized evolutions
to the canonical flow. We use $K,\psi,\Psi$ and $m$ from
Section~\ref{sec:flow} and fix $0<T<\infty$.

\subsection{Stability of the nonlinear resolvents}
\begin{lemma}\label{lem:reslimit}
For $g\in K$ and $\tau>0$, let $J_{\tau,s}$ and $J_{\tau,*}$
be the resolvents of Lemma~\ref{lem:resolvent}. Then
\begin{equation}\label{eq:reslim}
 J_{\tau,s}g\longrightarrow J_{\tau,*}g
 \quad\text{in }L^2(\rho),\qquad s\downarrow0.
\end{equation}
\end{lemma}
\begin{proof}
Put $v_s=(J_{\tau,s}g)^2$ and $v_*=(J_{\tau,*}g)^2$.
By \eqref{eq:J}, they minimize
\begin{equation}\label{eq:minimizers}
 F_s(v)=\int_E(\Psi(v)-gv)\,d\rho+\frac{\tau}{4}\E_s(v,v)
\end{equation}
on $\F_s$ and $\F_*$, respectively, with $s=*$ in the latter
case. Moreover, $c^2\le v_s,v_*\le C^2$.
Since $v_*\in\F_*\subset\F_s$ and
$\E_s(v_*,v_*)\uparrow\E_*(v_*,v_*)$, comparison with $v_*$ gives
\begin{equation}\label{eq:up}
 \limsup_{s\downarrow0}F_s(v_s)
 \le\lim_{s\downarrow0}F_s(v_*)=F_*(v_*).
\end{equation}
The density bounds give weak compactness in $L^2(\rho)$ and keep
the potential term in \eqref{eq:minimizers} bounded. Thus
\eqref{eq:up} also bounds $\E_s(v_s,v_s)$.
Along any weakly convergent subsequence $v_s\rightharpoonup v$,
Mosco lower semicontinuity puts $v\in\F_*$. Together with weak
lower semicontinuity of the convex potential, it gives
\begin{equation}\label{eq:lower}
 F_*(v)\le\liminf_{s\downarrow0}F_s(v_s)
 \le\limsup_{s\downarrow0}F_s(v_s)\le F_*(v_*).
\end{equation}
Strict convexity forces every such limit to equal $v_*$, and
the minimum values converge to $F_*(v_*)$.

To obtain strong convergence, use $v_*-v_s\in\F_s$ in the Euler
equation \eqref{eq:euler} for the form $\E_s$. The first variation
cancels, leaving
\begin{equation}\label{eq:strongmosco}
 \begin{aligned}
 F_s(v_*)-F_s(v_s)
 &=\int_E\bigl[\Psi(v_*)-\Psi(v_s)
                  -\psi(v_s)(v_*-v_s)\bigr]d\rho\\
 &\quad+\frac\tau4\E_s(v_*-v_s,v_*-v_s)\\
 &\ge\frac m2\|v_s-v_*\|_2^2.
 \end{aligned}
\end{equation}
By \eqref{eq:up} and \eqref{eq:lower}, the left side tends to zero.
Finally, the lower bound on the densities gives
$\|J_{\tau,s}g-J_{\tau,*}g\|_2
\le(2c)^{-1}\|v_s-v_*\|_2$, proving \eqref{eq:reslim}.
\end{proof}

\subsection{Stability of Gaussian gradients}
The next estimate identifies the weak Gaussian gradient from
convergence of the spectral functions. Convergence of the energy
norms then gives strong convergence.
\begin{lemma}\label{lem:gradient-stability}
Let $z_s\in L^2(0,T;\F_s)$ and $z\in L^2(0,T;\F_*)$. If
\begin{equation}\label{eq:gradient-stability-hypotheses}
 \begin{aligned}
 z_s&\longrightarrow z &&\text{in }L^2(0,T;L^2(\rho)),\\
 \int_0^T\E_s(z_{s,t},z_{s,t})\,dt
 &\longrightarrow\int_0^T\E_*(z_t,z_t)\,dt,
 \end{aligned}
\end{equation}
then
\begin{equation}\label{eq:gradient-stability}
 D_s z_s\longrightarrow\D(z(Y))
 \quad\text{in }L^2((0,T)\times\Omega,dt\,d\PP;H).
\end{equation}
\end{lemma}
\begin{proof}
The energy assumption and \eqref{eq:static-gradients} bound
$D_s z_s$ in $L^2(dt\,d\PP;H)$.
Let $\xi_k$ be the Gaussian coordinate along $e_k$, and let
$\partial_k$ denote differentiation in this direction.
For $\Phi=\varphi(\xi_1,\ldots,\xi_n)$ with
$\varphi\in C_c^\infty(\R^n)$ and
$\chi\in C_c^\infty(0,T)$, the weak coordinate derivative of
Lemma~\ref{lem:weighted-gradient-domain} gives
\begin{equation}\label{eq:coordinate-gradient-test}
 \begin{aligned}
 &\int_0^T\chi(t)\EE[\Phi\langle D_sz_{s,t},e_k\rangle_H]\,dt\\
 &\quad=q_k(s)^{1/2}\int_0^T
       \chi(t)\EE[z_{s,t}(Y)(\xi_k\Phi-\partial_k\Phi)]\,dt\\
 &\quad\longrightarrow
       \int_0^T\chi(t)\EE[z_t(Y)(\xi_k\Phi-\partial_k\Phi)]\,dt\\
 &\quad=\int_0^T\chi(t)
       \EE[\Phi\langle\D(z_t(Y)),e_k\rangle_H]\,dt.
 \end{aligned}
\end{equation}
Here $q_k(s)\to1$, and the lift by $Y$ is an isometry:
$\|z_s(Y)-z(Y)\|_{L^2(dt\,d\PP)}
 =\|z_s-z\|_{L^2(dt\,d\rho)}\to0$.
All the gradients are measurable with respect to the Gaussian field.
Finite sums of the tests $\chi\Phi e_k$ are dense in this Gaussian
$L^2$ subspace, so \eqref{eq:coordinate-gradient-test} and the
uniform bound imply
\[
 D_s z_s\rightharpoonup\D(z(Y))
 \quad\text{in }L^2(dt\,d\PP;H).
\]
The second condition in \eqref{eq:gradient-stability-hypotheses}
and \eqref{eq:static-canonical-gradient} give convergence of the
squared norms. Taking norms and inner products in
$L^2(dt\,d\PP;H)$, we conclude that
\[
 \begin{aligned}
 \|D_s z_s-\D(z(Y))\|_2^2
 &=\|D_s z_s\|_2^2+\|\D(z(Y))\|_2^2\\
 &\quad-2\langle D_s z_s,\D(z(Y))\rangle_{L^2(dt\,d\PP;H)}
 \longrightarrow0.
 \end{aligned}
\]
\end{proof}

\subsection{Conditional response estimates}
We now transfer Gaussian gradient estimates to the conditional
responses that enter the drift.
\begin{lemma}\label{lem:response-estimate}
For $g_s\in L^2(0,T;\F_s)$, $g\in L^2(0,T;\F_*)$ and $1\le i\le N$,
\begin{equation}\label{eq:responsebound}
 \begin{aligned}
 &\int_0^T\|R_{i,s}[g_{s,t}]-R_i[g_t]\|_{L^1(\rho)}\,dt\\
 &\quad\le\sqrt T\,\|D_s y_i\|_{L^2(\PP;H)}
      \|D_sg_s-\D(g(Y))\|_{L^2(dt\,d\PP;H)}\\
 &\qquad+\sqrt T\,\|D_s y_i-\D Y_i\|_{L^2(\PP;H)}
      \|\D(g(Y))\|_{L^2(dt\,d\PP;H)}.
 \end{aligned}
\end{equation}
For each $g\in\F_*$,
\begin{equation}\label{eq:joint-response}
 \|R_i[g]\|_{L^1(\rho)}
 \le\|\D Y_i\|_{L^2(\PP;H)}\,\E_*(g,g)^{1/2}.
\end{equation}
\end{lemma}
\begin{proof}
Subtract the conditional expectations defining the two responses.
Add and subtract the scalar product of $D_s y_i$ with
$\D(g_t(Y))$. Conditional Jensen and the triangle inequality give,
for almost every $t$,
\[
 \begin{aligned}
 \|R_{i,s}[g_{s,t}]-R_i[g_t]\|_{L^1(\rho)}
 &\le\EE\bigl[\|D_s y_i\|_H
                    \|D_sg_{s,t}-\D(g_t(Y))\|_H\bigr]\\
 &\quad+\EE\bigl[\|D_s y_i-\D Y_i\|_H
                         \|\D(g_t(Y))\|_H\bigr].
 \end{aligned}
\]
Integrate in time and use Cauchy--Schwarz on $(0,T)\times\Omega$.
Since $D_s y_i$ and $\D Y_i$ do not depend on time, this gives the
factor $\sqrt T$ in \eqref{eq:responsebound}.
The bound \eqref{eq:joint-response} follows from the
Cauchy--Schwarz bound after \eqref{eq:gamma} with $u=y_i$, since
$\E_*(y_i,y_i)=\EE\|\D Y_i\|_H^2$.
\end{proof}

For the remaining coefficient terms, let $u_s\to u$ in
$dt\,d\rho$-measure with $|u_s|,|u|\le B<\infty$.
For any $w\in L^1((0,T)\times E,dt\,d\rho)$,
\begin{equation}\label{eq:multiplier}
 \int_0^T\int_E|u_s-u|\,|w(t,y)|\,d\rho\,dt
 \longrightarrow0.
\end{equation}
Indeed, splitting the integral at $|w|=R$ bounds it by
\[
 R\|u_s-u\|_{L^1(dt\,d\rho)}
 +2B\int_0^T\int_{\{|w(t,\cdot)|>R\}}
                 |w(t,y)|\,d\rho(y)\,dt.
\]
The first term tends to zero as $s\downarrow0$ by bounded
convergence in measure on the finite measure space.
Then $R\uparrow\infty$ removes the integrable tail and proves
\eqref{eq:multiplier}.

\section{McKean--Vlasov weak realization}
\label{sec:realize}
To realize the canonical flow, we need integrability of the
law-dependent coefficients under that same law. The following
estimate supplies it from the response energy.

\subsection{An energy estimate for the coefficients}
\begin{lemma}\label{lem:coefficient-estimate}
Let $0<T<\infty$ and $g\in L^2(0,T;\F_*)$ with $c\le g\le C$. The fields
$a[g_t]$ and $b[g_t]$ have jointly Borel representatives on
$(0,T)\times E$, and
\begin{equation}\label{eq:coefficient-energy-estimate}
 \begin{aligned}
 &\int_0^T\int_E
       \bigl(|b[g_t]|+\Tr a[g_t]\bigr)g_t\,d\rho\,dt\\
 &\quad\le TC^2\bigl(\|b\|_{L^1(\rho)}
                         +\|\Tr a\|_{L^1(\rho)}\bigr)\\
 &\qquad+C\sqrt T\sum_{i=1}^N\|\D Y_i\|_{L^2(\PP;H)}
       \left(\int_0^T\E_*(g_t,g_t)\,dt\right)^{1/2}.
 \end{aligned}
\end{equation}
\end{lemma}
\begin{proof}
The map $R_i:\F_*\to L^1(\rho)$ is continuous by
\eqref{eq:joint-response}. Thus $t\mapsto g_t$ and
$t\mapsto R_i[g_t]$ are strongly measurable in $L^1(\rho)$ and
admit jointly Borel representatives, with $c\le g\le C$. Formula
\eqref{eq:law-coefficients} gives the corresponding representatives
of the coefficients. For almost every $t$, the same formula and
\eqref{eq:joint-response} yield
\[
 \begin{aligned}
 \int_E\bigl(|b[g_t]|+\Tr a[g_t]\bigr)g_t\,d\rho
 &\le\int_E g_t^2\bigl(|b|+\Tr a\bigr)\,d\rho
          +\sum_{i=1}^N\int_E g_t|R_i[g_t]|\,d\rho\\
 &\le C^2\bigl(\|b\|_1+\|\Tr a\|_1\bigr)
       +C\sum_{i=1}^N\|\D Y_i\|_{L^2(\PP;H)}
                          \E_*(g_t,g_t)^{1/2}.
 \end{aligned}
\]
Integration in time and Cauchy--Schwarz give
\eqref{eq:coefficient-energy-estimate}.
\end{proof}

\subsection{Final proof of Theorem~\ref{thm:limit}}
\label{subsec:final-limit}
Fix $0<T<\infty$. For $\tau>0$, write
$g_j=J_{\tau,*}^j f_0\in K$. Contractivity gives a comparison
of the Euler iterates:
\begin{equation}\label{eq:euler-comparison}
 \begin{aligned}
 \|J_{\tau,s}^k f_0-J_{\tau,*}^k f_0\|_1
 &\le\|J_{\tau,s}^{k-1}f_0-J_{\tau,*}^{k-1}f_0\|_1
       +\|(J_{\tau,s}-J_{\tau,*})g_{k-1}\|_1\\
 &\le\sum_{j=0}^{k-1}
       \|(J_{\tau,s}-J_{\tau,*})g_j\|_1
 \longrightarrow0
 \end{aligned}
\end{equation}
for each fixed $k$, by Lemma~\ref{lem:reslimit}.
For the continuous flows, the nonlinear operators of
Section~\ref{sec:flow} are accretive, their domains have the
common $L^1$ closure $K$, and their resolvent ranges contain $K$.
These are \eqref{eq:accretivity} and \eqref{eq:nonlinear-range}
for each form. Since Lemma~\ref{lem:reslimit} gives $L^1(\rho)$
resolvent convergence for every datum in $K$, the semigroup convergence
theorem \cite[Theorem 3.1]{BP} applies and yields
\begin{equation}\label{eq:flowlim}
 \sup_{0\le t\le T}\|f_{s,t}-f_t\|_{L^1(\rho)}
 \longrightarrow0.
\end{equation}
This proves \eqref{eq:thmlimit1}.

To identify the response drift, we first compare the gradients of
the squared densities. Put $v_{s,t}=f_{s,t}^2$ and $v_t=f_t^2$.
The uniform density bounds and \eqref{eq:flowlim} give
\begin{equation}\label{eq:squaresstrong}
 f_s\longrightarrow f,\qquad v_s\longrightarrow v
 \quad\text{in }L^2(0,T;L^2(\rho)).
\end{equation}
The energy identity \eqref{eq:thmflow3} turns the density error at
time $T$ into an energy norm error:
\begin{equation}\label{eq:energy-limit}
 \begin{aligned}
 &\left|\int_0^T\E_s(v_{s,t},v_{s,t})\,dt
              -\int_0^T\E_*(v_t,v_t)\,dt\right|\\
 &\quad=\frac23\left|\int_E(f_{s,T}^3-f_T^3)\,d\rho\right|\\
 &\quad\le2C^2\|f_{s,T}-f_T\|_1\longrightarrow0.
 \end{aligned}
\end{equation}
Equations \eqref{eq:squaresstrong} and \eqref{eq:energy-limit}
verify the hypotheses of Lemma~\ref{lem:gradient-stability} with
$z_s=v_s$ and $z=v$. Hence
\begin{equation}\label{eq:grad1}
 D_s v_s\longrightarrow\D(v(Y))
 \quad\text{in }L^2(dt\,d\PP;H).
\end{equation}

The lower density bound now allows us to return from $v_s$ to
$f_s$. For every $s>0$ and for $s=*$, the square-root chain rule gives
\begin{equation}\label{eq:density-gradient-bound}
 \int_0^T\E_s(f_{s,t},f_{s,t})\,dt
 \le\frac1{4c^2}\int_0^T
           \E_s(f_{s,t}^2,f_{s,t}^2)\,dt<\infty.
\end{equation}
Thus the density gradients are defined, and
\begin{equation}\label{eq:divchain}
 D_sf_s=\frac1{2f_s(Y)}D_sv_s,\qquad
 \D(f(Y))=\frac1{2f(Y)}\D(v(Y)).
\end{equation}
The inverse densities are bounded and converge in
$dt\,d\rho$-measure. Moreover, $\Gamma(v_t,v_t)\in L^1(dt\,d\rho)$.
Conditional expectation and \eqref{eq:multiplier} therefore give,
with norms in $L^2(dt\,d\PP;H)$,
\begin{equation}\label{eq:inverse-density-multiplier}
 \begin{aligned}
 &\left\|\left(\frac1{2f_s(Y)}-\frac1{2f(Y)}\right)
                        \D(v(Y))\right\|_2^2\\
 &\quad=\int_0^T\int_E
       \left|\frac1{2f_{s,t}}-\frac1{2f_t}\right|^2
                    \Gamma(v_t,v_t)\,d\rho\,dt
 \longrightarrow0.
 \end{aligned}
\end{equation}
Returning to \eqref{eq:divchain}, we obtain
\begin{equation}\label{eq:grad2}
 \begin{aligned}
 \|D_sf_s-\D(f(Y))\|_2
 &\le\frac1{2c}\|D_sv_s-\D(v(Y))\|_2\\
 &\quad+\left\|\left(\frac1{2f_s(Y)}-\frac1{2f(Y)}\right)
                         \D(v(Y))\right\|_2
 \longrightarrow0.
 \end{aligned}
\end{equation}

Since $Q_s^{1/2}\to I$ strongly and $\D Y_i\in L^2(\PP;H)$,
\eqref{eq:static-gradients} gives
\begin{equation}\label{eq:coordinate-response-limit}
 D_s y_i=Q_s^{1/2}\D Y_i\longrightarrow\D Y_i
 \quad\text{in }L^2(\PP;H).
\end{equation}
We may apply that lemma with $g_s=f_s$ and $g=f$ by
\eqref{eq:density-gradient-bound}. Its right side tends to zero
by \eqref{eq:grad2} and \eqref{eq:coordinate-response-limit}, proving
\begin{equation}\label{eq:rconv}
 \int_0^T\|R_{i,s}[f_{s,t}]-R_i[f_t]\|_{L^1(\rho)}\,dt
 \longrightarrow0,\qquad 1\le i\le N.
\end{equation}
This is \eqref{eq:thmlimit2}. Substituting the coefficient
formulas \eqref{eq:law-coefficients} and
\eqref{eq:regularized-law-coefficients} also gives
\begin{equation}\label{eq:coeffconv}
 \begin{aligned}
 &\int_0^T\bigl(\|a_s[f_{s,t}]-a[f_t]\|_1
                    +\|b_s[f_{s,t}]-b[f_t]\|_1\bigr)\,dt\\
 &\quad\le CT\bigl(\|a_s-a\|_1+\|b_s-b\|_1\bigr)\\
 &\qquad+\int_0^T\int_E
            |f_{s,t}-f_t|\bigl(|a|+|b|\bigr)\,d\rho\,dt\\
 &\qquad+\sum_{i=1}^N\int_0^T
           \|R_{i,s}[f_{s,t}]-R_i[f_t]\|_1\,dt
 \longrightarrow0.
 \end{aligned}
\end{equation}
Here \eqref{eq:static-coefficients} controls the static coefficient
error, \eqref{eq:multiplier} applies with weight $|a|+|b|$, and
\eqref{eq:rconv} controls the response error.

We now construct the process with the canonical marginal laws
$\mu_t=f_t\rho$. The continuity in Theorem~\ref{thm:flow} gives
\begin{equation}\label{eq:mucont}
 \|\mu_t-\mu_r\|_{\mathrm{TV}}
 \le\|f_t-f_r\|_1\longrightarrow0\qquad(r\to t).
\end{equation}
In particular, $(\mu_t)$ is narrowly continuous.
The canonical case of \eqref{eq:density-gradient-bound} permits
Lemma~\ref{lem:coefficient-estimate} with $g=f$, giving jointly
Borel coefficients and
\begin{equation}\label{eq:trevint}
 \int_0^T\int_E
       \bigl(|b[f_t]|+\Tr a[f_t]\bigr)\,d\mu_t\,dt<\infty.
\end{equation}
Time-dependent coefficients are identified $dt\,d\rho$-almost
everywhere and evaluated through these jointly Borel versions.
For $u\in C_c^\infty(\R^N)$, take $g=f_r$ in the weighted
generator identity \eqref{eq:weightedgen}. The chain rule then gives
\begin{equation}\label{eq:fpcalc}
 \begin{aligned}
 \int_E\ell[f_r]u\,d\mu_r
 &=-\int_E f_r\Gamma(u,f_r)\,d\rho\\
 &=-\frac12\E_*(f_r^2,u)
 \end{aligned}
\end{equation}
for almost every $r$. Inserting this into \eqref{eq:nlweak} yields
\begin{equation}\label{eq:fpe}
 \int_Eu\,d\mu_t-\int_Eu\,d\mu_0
 =\int_0^t\int_E
       \bigl(a[f_r]:\nabla^2u+b[f_r]\cdot\nabla u\bigr)\,d\mu_r\,dr.
\end{equation}
Time integration and smooth approximation, justified by
\eqref{eq:trevint}, extend \eqref{eq:fpe} to test functions in
$C_c^{1,2}((0,T)\times\R^N)$.
Choose a symmetric nonnegative Borel version of $a$, and extend
$a[f_t]$ and $b[f_t]$ by zero outside $E$. For the covariance
$2a[f_t]=2f_ta$ and drift $b[f_t]$, equations
\eqref{eq:mucont}, \eqref{eq:trevint}, and \eqref{eq:fpe}
verify the hypotheses of Trevisan's superposition theorem
\cite[Theorem 2.5]{Trevisan}.

The theorem supplies a probability law $\mathbb Q$ on
$C([0,T];\R^N)$. Its coordinate process $X$ has
$\Law_{\mathbb Q}(X_t)=\mu_t$ for every $t$, and
\begin{equation}\label{eq:martingale-problem}
 u(X_t)-u(X_0)-\int_0^t\ell[f_r]u(X_r)\,dr
 \quad\text{is a }\mathbb Q\text{-martingale}
\end{equation}
for every $u\in C_c^\infty(\R^N)$, in the usual augmented natural
filtration. All marginals are supported on the closed set $E$.
Applying this at rational times and using continuity gives
$X_t\in E$ for all $t\in[0,T]$, almost surely.
The marginal identities also give
\begin{equation}\label{eq:path-coefficient-integrability}
 \begin{aligned}
 &\EE_{\mathbb Q}\int_0^T
       \bigl(|b[f_t](X_t)|+\Tr a[f_t](X_t)\bigr)\,dt\\
 &\quad=\int_0^T\int_E
       \bigl(|b[f_t]|+\Tr a[f_t]\bigr)\,d\mu_t\,dt<\infty.
 \end{aligned}
\end{equation}
Thus the drift integrals are finite almost surely. It remains
to identify the covariance of their martingale remainder.

For integers $R\ge1$, choose $\chi_R\in C_c^\infty(\R^N)$
equal to one on $B_R$ and zero outside $B_{2R}$, and put
\begin{equation}\label{eq:coordinate-cutoffs}
 \begin{aligned}
 u_i^R(y)&=y_i\chi_R(y),\qquad
 u_{ij}^R(y)=y_i y_j\chi_R(y),\\
 \tau_R&=\inf\{t\in[0,T]:|X_t|\ge R\},
 \qquad \inf\varnothing=\infty.
 \end{aligned}
\end{equation}
Each continuous path is bounded on $[0,T]$, so $\tau_R>T$ for
all sufficiently large $R$.
On $B_R$, the generator acts as
\begin{equation}\label{eq:coordinate-generators}
 \ell[f_t]u_i^R=b_i[f_t],\qquad
 \ell[f_t]u_{ij}^R
   =y_i b_j[f_t]+y_j b_i[f_t]+2f_ta_{ij}.
\end{equation}
Define
\begin{equation}\label{eq:coordinate-martingales}
 M_i(t)=X_i(t)-X_i(0)-\int_0^t b_i[f_r](X_r)\,dr.
\end{equation}
Applying \eqref{eq:martingale-problem} to $u_i^R$ and stopping
at $\tau_R$ shows that $M_i(\cdot\wedge\tau_R)$ is a martingale.
Hence $M_i$ is a continuous local martingale.
The product formula, up to this stopping time, reads
\[
 \begin{aligned}
 d(X_iX_j)
 &=X_i\,dM_j+X_j\,dM_i+d\langle M_i,M_j\rangle_t\\
 &\quad+\bigl(y_i b_j[f_t]+y_j b_i[f_t]\bigr)(X_t)\,dt.
 \end{aligned}
\]
Subtracting the martingale identity for $u_{ij}^R$ shows that
$\langle M_i,M_j\rangle_{t\wedge\tau_R}
 -2\int_0^{t\wedge\tau_R}f_r(X_r)a_{ij}(X_r)\,dr$
is a continuous local martingale of finite variation.
It starts at zero and hence vanishes:
\begin{equation}\label{eq:stopped-brackets}
 \langle M_i,M_j\rangle_{t\wedge\tau_R}
 =2\int_0^{t\wedge\tau_R}f_r(X_r)a_{ij}(X_r)\,dr.
\end{equation}
Letting $R\uparrow\infty$ in \eqref{eq:stopped-brackets} removes
the stopping.
With the predictable symmetric nonnegative matrix
$\Sigma_t=\sqrt{f_t(X_t)}\,\sigma(X_t)$, the bracket becomes
\begin{equation}\label{eq:brackets}
 d\langle M\rangle_t
 =2f_t(X_t)a(X_t)\,dt=\Sigma_t\Sigma_t^T\,dt.
\end{equation}
By \eqref{eq:path-coefficient-integrability},
$\EE_{\mathbb Q}\Tr\langle M\rangle_T<\infty$, so $M$ is
square-integrable. The continuous-martingale representation theorem
\cite[Section 3.4.A]{KS} therefore gives, on an extension of the
filtered probability space, an $N$-dimensional Brownian motion $W$
such that
\begin{equation}\label{eq:bm}
 M_t=\int_0^t\Sigma_r\,dW_r.
\end{equation}
Substituting this representation into
\eqref{eq:coordinate-martingales} and using
\eqref{eq:law-coefficients} gives
\begin{equation}\label{eq:finalsde}
 \begin{aligned}
 X_i(t)
 &=X_i(0)+\int_0^t
       \bigl(f_r b_i+R_i[f_r]\bigr)(X_r)\,dr\\
 &\quad+\sum_{j=1}^N\int_0^t
       \sqrt{f_r(X_r)}\,\sigma_{ij}(X_r)\,dW_j(r),\\
 \Law(X_t)&=f_t\rho,\qquad 0\le t\le T.
 \end{aligned}
\end{equation}
This proves \eqref{eq:thmsde}.

\clearpage
\appendix
\section{Static spectral-response results}\label{sec:appendix}
This appendix proves the static inputs in Proposition~\ref{prop:static}.
The Green-energy estimates lead to the Gaussian response forms.
We then identify the compensated drift and the canonical increasing limit.
In this appendix, $C$ denotes a positive deterministic constant whose
value may change between estimates.

Use the basis $(e_k)$ of \eqref{eq:Q}. The Gaussian coordinates
and derivatives are
\begin{equation}\label{eq:appendix-basis}
 h=\sum_{k\ge1}\xi_ke_k,\qquad
 \partial_kF=\langle\D F,e_k\rangle_H.
\end{equation}
The coordinates $(\xi_k)$ are independent standard normals, and
the field expansion is distributional. With $\delta$ denoting the
Gaussian divergence, write $\LL_s=-\delta Q_s\D$.
The random Green operator $T=T(h)$ is
\begin{equation}\label{eq:Tdef}
 \langle Tu,v\rangle_H=\int_Duv\,dM_h,\qquad
 Tv_i=\lambda_i^{-1}v_i,\qquad v_i=\lambda_i^{-1/2}\phi_i.
\end{equation}
The Gaussian Sobolev spaces are denoted by $\mathbb D^{m,p}$.
By \cite[Theorem 1.8]{Cai} and the logarithmic change of variables, write
$\rho(dy)=p_N(y)\,dy$ on the ordered chamber, with $p_N$ extended
by zero to $\R^N$.

\subsection{Annealed Green-energy bounds}
\label{app:annealed}
The Gaussian gradients and compensated drift require moments of
$I_G(\nu_i)$. We obtain them from a ball-mass bound for $M_h$
and the normalization $\nu_i(D)=1$.
By \cite[proof of Lemma 6.10, equation (6.33)]{Cai}, for every
$1<q<4/\gamma^2$ and every dyadic square $Q$ of side $2^{-k}$,
\begin{equation}\label{eq:dyadic-mass-moment}
 \EE M_h(Q\cap D)^q\le C_q2^{-k\zeta(q)},\qquad
 \zeta(q)=\left(2+\frac{\gamma^2}{2}\right)q
                         -\frac{\gamma^2}{2}q^2,\quad k\ge0,
\end{equation}
where $C_q$ is independent of the position of $Q$.
Since $\zeta(1)=2$ and $\zeta'(1)=2-\gamma^2/2>0$, we may fix
$q>1$ such that $q<4/\gamma^2$ and $\zeta(q)>2$, and then choose
\[
 0<\alpha<\frac{\zeta(q)-2}{q}.
\]
These parameters remain fixed throughout the appendix.

\begin{lemma}\label{lem:frostman}
There is a random $C_\alpha\ge1+M_h(D)$ such that
\begin{equation}\label{eq:Frostman}
 M_h(B(x,r)\cap D)\le C_\alpha r^\alpha,
 \qquad x\in D,\quad 0<r\le1,\qquad
 \EE C_\alpha^q<\infty.
\end{equation}
For every fixed $i$ and every $p>0$,
\begin{equation}\label{eq:logmoments}
 \EE\lambda_i^p+\EE|Y_i|^p<\infty.
\end{equation}
\end{lemma}
\begin{proof}
Let $\mathcal Q_k$ be the dyadic squares of side $2^{-k}$ meeting $D$,
and put
\[
 Z=\sup_{k\ge0}\max_{Q\in\mathcal Q_k}
                         2^{k\alpha}M_h(Q\cap D).
\]
There are at most $C_D2^{2k}$ squares in $\mathcal Q_k$.
Thus \eqref{eq:dyadic-mass-moment} gives
\[
 \EE Z^q
 \le\sum_{k\ge0}\sum_{Q\in\mathcal Q_k}
                         2^{k\alpha q}\EE M_h(Q\cap D)^q
 \le C\sum_{k\ge0}2^{-k(\zeta(q)-2-\alpha q)}<\infty.
\]
For $2^{-j-1}<r\le2^{-j}$, the ball $B(x,r)$ meets at most
$16$ dyadic squares of side $2^{-j}$. Consequently,
\[
 M_h(B(x,r)\cap D)\le16Z2^{-j\alpha}
                         \le16\cdot2^\alpha Zr^\alpha.
\]
Also $M_h(D)\le\#\mathcal Q_0\,Z$.
Taking $C_\alpha=1+M_h(D)+16\cdot2^\alpha Z$ proves
\eqref{eq:Frostman}.

To bound $\lambda_i$ from above, choose smooth functions
$g_1,\ldots,g_i$ with disjoint compact supports in $D$, each equal
to one on a ball $B_l\Subset D$. Their Dirichlet energies and their
$L^2(M_h)$ masses are diagonal on this span, so min--max yields
\[
 \lambda_i\le\max_{l\le i}
       \frac{\|g_l\|_H^2}{\int_Dg_l^2\,dM_h}
       \le C_i\max_{l\le i}M_h(B_l)^{-1}.
\]
The interior ball masses have all negative moments
\cite[Theorem 2.12]{RV}. Hence
\[
 \EE\lambda_i^p\le C_{i,p}\sum_{l\le i}
                            \EE M_h(B_l)^{-p}<\infty,
 \qquad p>0.
\]

For the logarithmic moments we also need a lower bound on
$\lambda_1$. The Green bound
$G_D(x,z)\le C_D+\log^+(1/|x-z|)$ and the layer-cake formula give,
for $m=1,2$,
\begin{equation}\label{eq:green-frostman-potential}
 \begin{aligned}
 K_m&:=\sup_{x\in D}\int_DG_D(x,z)^m\,M_h(dz)\\
 &\le C_{D,m}\left(M_h(D)+
       \int_0^1(1+\log(1/r))^{m-1}
           \sup_{x\in D}M_h(B(x,r)\cap D)\,\frac{dr}{r}\right)\\
 &\le C_{D,\alpha,m}C_\alpha.
 \end{aligned}
\end{equation}
Indeed, the last integral is bounded by
$C_\alpha\int_0^1(1+\log(1/r))^{m-1}r^{\alpha-1}\,dr<\infty$.
The Green operator on $L^2(M_h)$ has largest eigenvalue
$\lambda_1^{-1}$. Schur's estimate therefore gives
$\lambda_1^{-1}\le K_1\le C_{D,\alpha}C_\alpha$.
Since $\lambda_i\ge\lambda_1$,
\[
 |Y_i|\le\log^+\lambda_i+
                         \log^+(C_{D,\alpha}C_\alpha).
\]
For every $p>0$, it follows that
\[
 \EE|Y_i|^p\le C_{D,\alpha,p,q}
           \bigl(1+\EE\lambda_i^p+\EE C_\alpha^q\bigr)<\infty,
\]
which completes \eqref{eq:logmoments}.
\end{proof}

\begin{proposition}\label{prop:AG}
Fix $0<\theta<\min\{1,q/4\}$. For every fixed $i$ there is a
random $C_{\alpha,i}\ge1$ with $\EE C_{\alpha,i}^\theta<\infty$ such that
\begin{equation}\label{eq:spectralmass}
 \begin{aligned}
 \nu_i(B(x,r))&\le\min\{1,C_{\alpha,i} r^\alpha\},
                  \qquad x\in D,\quad 0<r\le1,\\
 I_G(\nu_i)&\le C_{D,\alpha}(1+\log C_{\alpha,i}).
 \end{aligned}
\end{equation}
Consequently $\EE I_G(\nu_i)^p<\infty$ for every $1\le p<\infty$.
\end{proposition}
\begin{proof}
Use the continuous representative of $\phi_i$ and its Green
representation from \cite[Lemma 6.27]{Cai}. Cauchy--Schwarz,
$\int_D\phi_i^2\,dM_h=1$, and
\eqref{eq:green-frostman-potential} give
\[
 \begin{aligned}
 |\phi_i(x)|^2
 &=\lambda_i^2\left|\int_DG_D(x,z)\phi_i(z)\,M_h(dz)\right|^2\\
 &\le\lambda_i^2\int_DG_D(x,z)^2\,M_h(dz)
 \le C_{D,\alpha}\lambda_i^2C_\alpha.
 \end{aligned}
\]
The ball-mass estimate \eqref{eq:Frostman} now gives
\[
 \nu_i(B(x,r))
 \le\|\phi_i\|_\infty^2M_h(B(x,r)\cap D)
 \le C_{D,\alpha}\lambda_i^2C_\alpha^2r^\alpha.
\]
Since $\nu_i$ is a probability measure, the first bound in
\eqref{eq:spectralmass} holds with
$C_{\alpha,i}=1+C_{D,\alpha}\lambda_i^2C_\alpha^2$.
By Cauchy--Schwarz,
\[
 \begin{aligned}
 \EE C_{\alpha,i}^\theta
 &\le C_\theta\bigl(1+
                       \EE[\lambda_i^{2\theta}C_\alpha^{2\theta}]\bigr)\\
 &\le C_\theta\bigl(1+
       (\EE\lambda_i^{4\theta})^{1/2}
       (\EE C_\alpha^{4\theta})^{1/2}\bigr)<\infty.
 \end{aligned}
\]
Both moments are finite by Lemma~\ref{lem:frostman}, since
$4\theta<q$.

To estimate the Green energy, apply the logarithmic layer-cake
formula to $\nu_i$. Its total mass is one, and its small-ball
masses are bounded by \eqref{eq:spectralmass}, so
\[
 \begin{aligned}
 \sup_{x\in D}\int_DG_D(x,z)\,\nu_i(dz)
 &\le C_D+\int_0^1\min\{1,C_{\alpha,i} r^\alpha\}\,\frac{dr}{r}\\
 &=C_D+\frac{1+\log C_{\alpha,i}}{\alpha}.
 \end{aligned}
\]
The equality follows by splitting the integral at
$r=C_{\alpha,i}^{-1/\alpha}$. Integrating once more against $\nu_i$
gives the Green-energy bound in \eqref{eq:spectralmass}.
Finally, for $1\le p<\infty$,
\[
 \EE I_G(\nu_i)^p
 \le C_{D,\alpha,p}\EE(1+\log C_{\alpha,i})^p
 \le C_{D,\alpha,p,\theta}\EE C_{\alpha,i}^\theta<\infty.
\]
\end{proof}

\begin{remark}\label{rem:near-diagonal}
For the exponent $\theta$ of Proposition~\ref{prop:AG} and
$0<r\le1$, the same ball-mass bound gives
\[
 \begin{aligned}
 \EE(\nu_i\otimes\nu_i)\{|x-z|<r\}
 &=\EE\int_D\nu_i(B(x,r))\,\nu_i(dx)\\
 &\le\EE\min\{1,C_{\alpha,i} r^\alpha\}
 \le r^{\alpha\theta}\EE C_{\alpha,i}^\theta.
 \end{aligned}
\]
\end{remark}

\subsection{Gaussian gradients and response forms}
\label{app:forms}
We combine the Green-energy moments of Proposition~\ref{prop:AG}
with the pathwise first variation of the eigenvalues. To carry out
Gaussian integration by parts, we first choose a version of the GMC
compatible with finite Cameron--Martin shifts. Put
$V=\bigcup_{M\ge1}\operatorname{span}\{e_1,\ldots,e_M\}$.

\begin{lemma}\label{lem:coherent-shifts}
There is a $V$-invariant Gaussian full-measure set $\Omega_*$ such
that, for $h\in\Omega_*$ and $g\in V$,
\begin{equation}\label{eq:alltilts}
 h+g\in\Omega_*,\qquad M_{h+g}=e^{\gamma g}M_h.
\end{equation}
Each $M_h$ is finite, has full support, and satisfies the
ball-mass bound in \eqref{eq:Frostman}.
\end{lemma}
\begin{proof}
Write $P_\varepsilon^D=e^{\varepsilon\Delta_D}$ for the heat
regularization and $M_{h,\varepsilon}$ for the full-variance
GMC approximants of \cite[(6.9)]{Cai}. Fix the deterministic
subsequence of \cite[Lemma 6.2]{Cai}. For $g\in V$, they satisfy
\begin{equation}\label{eq:coherent-regularized-shift}
 M_{h+g,\varepsilon_m}
 =e^{\gamma P_{\varepsilon_m}^Dg}M_{h,\varepsilon_m},
 \qquad
 \|P_{\varepsilon_m}^Dg-g\|_\infty\longrightarrow0.
\end{equation}
The uniform convergence holds because $g$ is a finite combination
of Dirichlet eigenfunctions. Whenever $M_{h,\varepsilon_m}$
converges, \eqref{eq:coherent-regularized-shift} forces
$M_{h+g,\varepsilon_m}$ to converge to $e^{\gamma g}M_h$.
Applying the same identity with $-g$ gives the converse. Thus the
convergence event is invariant under every shift in $V$.
On this event,
\[
 e^{-\gamma\|g\|_\infty}M_h\le M_{h+g}
                    \le e^{\gamma\|g\|_\infty}M_h.
\]
The comparison preserves finiteness, full support, and the
ball-mass bound. Intersecting the convergence event with the
full-measure event of these properties gives $\Omega_*$.
\end{proof}

For $g\in V$, min--max applied to \eqref{eq:alltilts} gives
\begin{equation}\label{eq:direction-lipschitz}
 |Y_i(h+\tau g)-Y_i(h+\tau' g)|
       \le\gamma\|g\|_\infty|\tau-\tau'|.
\end{equation}
The tilted mass form is analytic on $H$, since the multiplier
$e^{\gamma\tau g}$ is analytic in $L^\infty(D)$.
At a simple eigenvalue we may therefore choose a differentiable
branch $\lambda_i(\tau),\phi_i(\tau)$, with $\phi_i(\tau)$
normalized in $L^2(M_{h+\tau g})$.
Writing dots for derivatives at $\tau=0$, differentiation of the
normalization and \eqref{eq:spectrum} yields
\[
 \begin{aligned}
 0&=2\int_D\phi_i\dot\phi_i\,dM_h
                          +\gamma\int_Dg\,d\nu_i,\\
 \dot\lambda_i
 &=2\langle\phi_i,\dot\phi_i\rangle_H
   =2\lambda_i\int_D\phi_i\dot\phi_i\,dM_h.
 \end{aligned}
\]
Eliminating the eigenfunction derivative gives the first response
\cite[Corollary 6.35]{Cai}:
\begin{equation}\label{eq:firstresponse}
 \partial_gY_i=\frac{\dot\lambda_i}{\lambda_i}
                         =-\gamma\int_Dg\,d\nu_i.
\end{equation}
The Green--Riesz representation of this functional
\cite[Corollary A.5]{Cai}, together with the annealed energy
bounds, now identifies its closed Gaussian gradient.

\begin{lemma}\label{lem:MalliavinY}
For every $1<p<\infty$, $Y_i\in\mathbb D^{1,p}$ and
$\D Y_i=-\gamma U_D\nu_i$. The random response matrices satisfy
\begin{equation}\label{eq:covconvergence}
 0\le A_s\le A_0,\qquad
 A_{ij,s}\longrightarrow A_{ij,0}
       \quad\text{in }L^p(\PP),\quad 1\le p<\infty,
       \quad s\downarrow0.
\end{equation}
The matrix inequalities are in the sense of quadratic forms.
\end{lemma}
\begin{proof}
Fix $1<p<\infty$. Proposition~\ref{prop:AG} gives
\[
 \EE\|U_D\nu_i\|_H^p=\EE I_G(\nu_i)^{p/2}<\infty.
\]
By \eqref{eq:direction-lipschitz}, the restriction of $Y_i$ to
almost every $e_k$-coordinate line is absolutely continuous.
Almost-sure simplicity and Fubini's theorem then give, by
\eqref{eq:firstresponse},
\[
 \partial_kY_i=-\gamma\int_De_k\,d\nu_i
                   =-\gamma\langle U_D\nu_i,e_k\rangle_H\quad\text{a.e.}
\]
For every smooth compactly supported cylinder test $\Phi$,
Gaussian integration by parts on that line yields
\begin{equation}\label{eq:Gaussian-weak-gradient}
 \EE[Y_i(\xi_k\Phi-\partial_k\Phi)]
             =-\gamma\EE[\langle U_D\nu_i,e_k\rangle_H\Phi].
\end{equation}
These integrals are finite by the moment bounds for $Y_i$ and
$U_D\nu_i$.

To pass from the coordinate derivatives to the closed gradient,
let $\mathcal G_m=\sigma(\xi_1,\ldots,\xi_m)$ and
$Y_i^{(m)}=\EE[Y_i\mid\mathcal G_m]$. Denote by $\Pi_m$ the
orthogonal projection onto $\operatorname{span}\{e_1,\ldots,e_m\}$.
Testing \eqref{eq:Gaussian-weak-gradient} with
$\mathcal G_m$-measurable cylinders identifies the finite-dimensional
weak gradient:
\begin{equation}\label{eq:conditional-gradient}
 \D Y_i^{(m)}=-\gamma\Pi_m\EE[U_D\nu_i\mid\mathcal G_m].
\end{equation}
Finite-dimensional Gaussian Sobolev theory gives
$Y_i^{(m)}\in\mathbb D^{1,p}$. Since
$\mathcal G_m\uparrow\sigma(h)$, martingale convergence gives
$Y_i^{(m)}\to Y_i$ in $L^p(\PP)$. By
\eqref{eq:conditional-gradient},
\[
 \begin{aligned}
 \|\D Y_i^{(m)}+\gamma U_D\nu_i\|_{L^p(\PP;H)}
 &\le\gamma\|\EE[U_D\nu_i\mid\mathcal G_m]-U_D\nu_i\|_{L^p(\PP;H)}\\
 &\quad+\gamma\|(\Pi_m-I)U_D\nu_i\|_{L^p(\PP;H)}
 \longrightarrow0.
 \end{aligned}
\]
Closedness of $\D$ proves $Y_i\in\mathbb D^{1,p}$ with
$\D Y_i=-\gamma U_D\nu_i$.

For $z\in\R^N$,
\[
 0\le z^TA_sz
 =\left\|Q_s^{1/2}\sum_i z_i\D Y_i\right\|_H^2
 \le\left\|\sum_i z_i\D Y_i\right\|_H^2=z^TA_0z.
\]
Also $Q_s^{1/2}\to I$ strongly on $H$, and
\[
 |A_{ij,s}|\le\|\D Y_i\|_H\|\D Y_j\|_H,\qquad
 \EE(\|\D Y_i\|_H\|\D Y_j\|_H)^p<\infty
 \quad(1\le p<\infty).
\]
Dominated convergence therefore proves the entrywise limits in
\eqref{eq:covconvergence}.
\end{proof}

We next close the weighted response operator on the spectral
space $L^2(\rho)$.
\begin{proposition}\label{prop:forms}
For every $s\ge0$, the preform \eqref{eq:form_s} is closable.
Its minimal closure $(\E_s,\F_s)$ is a regular, strongly local,
conservative, symmetric Dirichlet form on $E$.
Writing $D_s$ for the closure of
$u\mapsto Q_s^{1/2}\D(u(Y))$ on $\Ccore$, one has
\begin{equation}\label{eq:appendix-Ds}
 \E_s(u,v)=\EE\langle D_su,D_sv\rangle_H,
 \qquad u,v\in\F_s.
\end{equation}
Moreover, $1,y_i\in\F_s$ and
\begin{equation}\label{eq:appendix-coordinate-gradient}
 D_s1=0,\qquad D_sy_i=Q_s^{1/2}\D Y_i,
 \qquad 1\le i\le N.
\end{equation}
\end{proposition}
\begin{proof}
Take $u_\ell\in\Ccore$ with $u_\ell\to0$ in $L^2(\rho)$ and
$Q_s^{1/2}\D(u_\ell(Y))\to G$ in $L^2(\PP;H)$.
For every fixed $k$, diagonalization of $Q_s$ gives
\[
 \partial_k(u_\ell(Y))
       \longrightarrow q_k(s)^{-1/2}\langle G,e_k\rangle_H
       \quad\text{in }L^2(\PP).
\]
Since $u_\ell(Y)\to0$ in $L^2(\PP)$ and $\partial_k$ is closed,
$\langle G,e_k\rangle_H=0$. Completeness of $(e_k)$ gives $G=0$.
Thus the response operator is closable; its closure $D_s$ has
\eqref{eq:appendix-Ds} by polarization.

For a smooth contraction $\eta$ with $\eta(0)=0$, the chain rule
on $\Ccore$ gives
\[
 Q_s^{1/2}\D((\eta\circ u)(Y))
                  =\eta'(u(Y))Q_s^{1/2}\D(u(Y)),\qquad
 \E_s(\eta\circ u,\eta\circ u)\le\E_s(u,u).
\]
Graph approximation and smooth approximation of normal contractions
therefore give the Markov property of the closure.
Symmetry and strong locality pass from the gradient preform to
its closure. Since $\Ccore$ is uniformly dense in $C_0(E)$ and
is a form core by construction, the closure is regular.

It remains to identify the constants and coordinate gradients.
Choose $0\le\chi_R\le1$ smooth, equal to one on $B_R$ and zero
outside $B_{2R}$, with $|\nabla\chi_R|\le C/R$.
Then $\chi_R\to1$ in $L^2(\rho)$ and
\[
 \|D_s\chi_R\|_{L^2(\PP;H)}^2
 =\E_s(\chi_R,\chi_R)
 \le\frac{C}{R^2}\EE\Tr A_0\longrightarrow0.
\]
Closedness gives $1\in\F_s$ and $D_s1=0$.

For the coordinate cutoff, the product rule on $\Ccore$ gives
\[
 \begin{aligned}
 D_s(y_i\chi_R)-Q_s^{1/2}\D Y_i
 &=(\chi_R(Y)-1)Q_s^{1/2}\D Y_i\\
 &\quad+Y_iQ_s^{1/2}\D(\chi_R(Y)).
 \end{aligned}
\]
Whenever $\nabla\chi_R(Y)\ne0$, one has $R\le|Y|\le2R$.
The chain rule, Cauchy--Schwarz, and $Q_s\le I$ therefore give
\begin{equation}\label{eq:forms-cutoff-bounds}
 \begin{aligned}
 \|D_s(y_i\chi_R)-Q_s^{1/2}\D Y_i\|_{L^2(\PP;H)}^2
 &\le2\EE[(1-\chi_R(Y))^2\|\D Y_i\|_H^2]\\
 &\quad+2\EE[Y_i^2|\nabla\chi_R(Y)|^2\Tr A_0]\\
 &\le C\EE[\one_{\{|Y|>R\}}\Tr A_0]\longrightarrow0.
 \end{aligned}
\end{equation}
Here $\Tr A_0\in L^1(\PP)$ by Proposition~\ref{prop:AG} and
Lemma~\ref{lem:MalliavinY}. Equation~\eqref{eq:forms-cutoff-bounds}
and $y_i\chi_R\to y_i$ in $L^2(\rho)$ allow us to apply
closedness of $D_s$, proving the coordinate assertion
in \eqref{eq:appendix-coordinate-gradient}.
Finally, $1\in\F_s$ and $\E_s(1,1)=0$ imply conservativity.
\end{proof}

At $s=0$ this construction gives the minimal closure
$(\E_0^{\min},\F_0^{\min})$; in statements indexed by $s\ge0$,
$(\E_0,\F_0)$ denotes this same form.

\subsection{Second spectral response and coordinate integration by parts}
\label{app:curvature}
Put $S_n=\sum_{i=1}^nY_i$ for $n\ge1$. For $1\le n\le N$, the
weak drift of the response form is obtained by integrating each term in
\[
 \E_s\Bigl(\sum_{i\le n}y_i,u\Bigr)
 =\sum_{k\ge1}q_k(s)\EE[\partial_kS_n\,\partial_k(u(Y))]
\]
against the corresponding Gaussian coordinate. Write
$h^{\perp k}=h-\xi_ke_k$ and let $\PP^{\perp k}$ be its law.
We will prove that, for every $k$ and $i$,
\begin{equation}\label{eq:linewise-target}
 \PP^{\perp k}\left\{h^{\perp k}:
  \exists t\in\R,\ 
  \lambda_i(h^{\perp k}+te_k)
   =\lambda_{i+1}(h^{\perp k}+te_k)\right\}=0.
\end{equation}
Along these collision-free lines, analytic perturbation gives the
classical second response. We first establish analyticity on finite
Gaussian fibers, then use two independent response directions to
bound the dimension of the collision set.

\subsubsection{Hilbert--Schmidt regularity on coherent fibers}
For a finite measure $M_h$ satisfying \eqref{eq:Frostman}, let
$J$ initially denote restriction from $C_c^\infty(D)$ to
$L^2(M_h)$. The estimate below extends $J$ to $H$, using
quasi-continuous representatives, and identifies $T(h)=J^*J$
for $h\in\Omega_*$ in \eqref{eq:Tdef}.

\begin{lemma}\label{lem:coherent}
For every $h\in\Omega_*$ of Lemma~\ref{lem:coherent-shifts},
$T(h)$ is compact and Hilbert--Schmidt, with
\begin{equation}\label{eq:coherent-HS}
 \|T(h)\|_{\mathrm{HS}}^2
 =\iint_{D^2}G_D(x,z)^2\,M_h(dx)M_h(dz)<\infty.
\end{equation}
For each $M\ge1$, the map
\begin{equation}\label{eq:analytic-operator-fiber}
 (t_1,\ldots,t_M)\longmapsto
 T\left(h+\sum_{l=1}^M t_le_l\right)
\end{equation}
is real analytic with values in $\mathrm{HS}(H)$.
\end{lemma}
\begin{proof}
The quantities $K_1,K_2$ in \eqref{eq:green-frostman-potential}
satisfy
\begin{equation}\label{eq:coherent-green-square}
 K_1+K_2\le C_{D,\alpha}C_\alpha.
\end{equation}
Green duality and Schur's estimate therefore give
\[
 \left|\int_Duv\,dM_h\right|
 \le\|u\|_H I_G(vM_h)^{1/2}
 \le K_1^{1/2}\|u\|_H\|v\|_{L^2(M_h)},
 \qquad u\in C_c^\infty(D),\quad v\in L^2(M_h).
\]
Taking the supremum over $v$ extends $J$ boundedly to $H$.
The operator $JJ^*$ on $L^2(M_h)$ has kernel $G_D(x,z)$,
which belongs to $L^2(M_h\otimes M_h)$ by
\eqref{eq:coherent-green-square}. Thus $JJ^*$ is Hilbert--Schmidt.
Its nonzero eigenvalues, with multiplicities, are those of
$J^*J=T(h)$, so
\begin{equation}\label{eq:T-HS-coherent}
 \|T(h)\|_{\mathrm{HS}}^2
 =\|JJ^*\|_{\mathrm{HS}}^2
 =\iint_{D^2}G_D(x,z)^2\,M_h(dx)M_h(dz)
 \le K_2M_h(D)<\infty.
\end{equation}

For a bounded real function $g$, put $T_{gM_h}u=J^*(gJu)$.
The same kernel
calculation yields
\[
 \begin{aligned}
 \|T_{gM_h}\|_{\mathrm{HS}}^2
 &=\iint_{D^2}g(x)g(z)G_D(x,z)^2\,M_h(dx)M_h(dz)\\
 &\le\|g\|_\infty^2\|T(h)\|_{\mathrm{HS}}^2.
 \end{aligned}
\]
Now fix $M$ and set $g_t=\sum_{l\le M}t_le_l$. The coherent-shift
identity \eqref{eq:alltilts} gives
\begin{equation}\label{eq:coherent-operator-series}
 \begin{aligned}
 T(h+g_t)&=\sum_{r=0}^\infty\frac{\gamma^r}{r!}
                                      T_{g_t^rM_h},\\
 \sum_{r=0}^\infty\frac{\gamma^r}{r!}
       \|T_{g_t^rM_h}\|_{\mathrm{HS}}
 &\le e^{\gamma\|g_t\|_\infty}\|T(h)\|_{\mathrm{HS}}.
 \end{aligned}
\end{equation}
Since $\|g_t\|_\infty$ is bounded on compact parameter sets,
this series converges locally uniformly in $\mathrm{HS}(H)$.
Each term is a homogeneous polynomial in $t$, proving
\eqref{eq:analytic-operator-fiber}.
\end{proof}

\subsubsection{Weak unique continuation and two response directions}
\label{subsec:unique-continuation}
Two independent responses will follow from the absence of quadratic
relations between orthonormal eigenfunctions. We first prove the
unique-continuation statement used to exclude such relations.

\begin{lemma}\label{lem:ucp}
Let $M_h$ satisfy \eqref{eq:Frostman}. If $\lambda>0$ and a locally
bounded continuous $u\in H^1_{\mathrm{loc}}(D)$ satisfies
\begin{equation}\label{eq:measure-eigenfunction}
 -\frac1{2\pi}\Delta u=\lambda uM_h
 \qquad\text{in distributions on }D,
\end{equation}
then $u$ vanishes identically whenever it vanishes on a nonempty
open subset of $D$.
\end{lemma}
\begin{proof}
Fix a disk $B=B(x_0,r)\Subset D$ with $r\le1/2$, and let $G_B$
be its Green function with the normalization of \eqref{eq:H}.
For bounded $f$ define
\[
 K_Bf(x)=\int_BG_B(x,z)f(z)\,M_h(dz),
 \qquad
 M_B=\sup_{x\in B}\int_BG_B(x,z)\,M_h(dz).
\]
Since $0\le G_B(x,z)\le C\log(2r/|x-z|)$, the layer-cake
formula gives
\begin{equation}\label{eq:local-green-schur}
 \begin{aligned}
 M_B
 &\le C\sup_{x\in B}
           \int_0^{2r}\frac{M_h(B(x,t)\cap D)}{t}\,dt\\
 &\le CC_\alpha\int_0^{2r}t^{\alpha-1}\,dt
 \le C_{D,\alpha}C_\alpha r^\alpha.
 \end{aligned}
\end{equation}
We may therefore choose the disks uniformly small enough that
\begin{equation}\label{eq:gauge-smallness}
 \lambda M_B<\frac12.
\end{equation}
Schur's estimate gives
$\|K_B\|_{L^2(M_h|_B)\to L^2(M_h|_B)}\le M_B$.
Green duality, as in the preceding proof, then yields the trace bound
\begin{equation}\label{eq:local-trace-bound}
 \int_B\psi^2\,dM_h
 \le\frac{M_B}{2\pi}\int_B|\nabla\psi|^2\,dx,
 \qquad \psi\in H_0^1(B),
\end{equation}
first for smooth $\psi$ and then by density.

On $L^\infty(M_h|_B)$ the operator $\lambda K_B$ has norm
less than $1/2$. Its Neumann series defines a positive gauge
\begin{equation}\label{eq:gauge-definition}
 w=\sum_{j=0}^\infty(\lambda K_B)^j1,
 \qquad w=1+\lambda K_Bw,
 \qquad 1\le w\le(1-\lambda M_B)^{-1}\le2.
\end{equation}
The trace bound implies that $wM_h|_B$ has finite $G_B$-energy.
Consequently $w-1=\lambda K_Bw\in H_0^1(B)$ and
\begin{equation}\label{eq:gauge-weak-equation}
 \frac1{2\pi}\int_B\nabla w\cdot\nabla\psi\,dx
   =\lambda\int_Bw\psi\,dM_h,
 \qquad \psi\in H_0^1(B).
\end{equation}
The potential in \eqref{eq:gauge-definition} is continuous. Indeed,
its near-diagonal contribution satisfies
\[
 \sup_{x\in B}\int_{B(x,\varepsilon)\cap B}
       G_B(x,z)|w(z)|\,M_h(dz)
 \le CC_\alpha\varepsilon^\alpha(1+|\log\varepsilon|)
 \longrightarrow0.
\]
Fix a small ball around a point of $B$. This estimate controls the
integral over that ball uniformly at nearby points; on its fixed
complement, continuity follows by dominated convergence.
Thus $w$ has a continuous representative satisfying $1\le w\le2$.

Set $v=u/w$. The product and reciprocal rules give
\begin{equation}\label{eq:gauge-quotient}
 v\in H^1_{\mathrm{loc}}(B),\qquad
 \nabla v=\frac{w\nabla u-u\nabla w}{w^2}
 \quad\text{a.e. on }B.
\end{equation}
For $\psi\in C_c^\infty(B)$, the functions $w\psi$ and $u\psi$
belong to $H_0^1(B)$. They are admissible tests in
\eqref{eq:measure-eigenfunction} and \eqref{eq:gauge-weak-equation}
by \eqref{eq:local-trace-bound}. Subtracting the two equations
eliminates the measure term:
\begin{equation}\label{eq:gauge-identity}
 \begin{aligned}
 0&=\int_B\nabla u\cdot\nabla(w\psi)\,dx
       -\int_B\nabla w\cdot\nabla(u\psi)\,dx\\
  &=\int_B(w\nabla u-u\nabla w)\cdot\nabla\psi\,dx\\
  &=\int_Bw^2\nabla v\cdot\nabla\psi\,dx.
 \end{aligned}
\end{equation}
Hence $v$ solves
\begin{equation}\label{eq:gauge-divergence}
 -\operatorname{div}(A\nabla v)=0,
 \qquad A=w^2I,\qquad I\le A\le4I.
\end{equation}
The matrix $A$ is bounded, measurable, symmetric and uniformly
elliptic. The planar unique-continuation result
\cite[Lemma 2.3]{HPS}, with zero lower-order term, applies to
$v\in H^1_{\mathrm{loc}}(B)$. If $u$ vanishes on a nonempty open
subset of $B$, so does $v$, and hence $v=0$ throughout $B$.
Since $w>0$, we obtain $u=0$ on $B$. Overlapping small disks
propagate the vanishing throughout the connected domain $D$.
\end{proof}

Fix a positive eigenvalue $\lambda$ of \eqref{eq:spectrum} with real
eigenspace $V_\lambda\subset L^2(M_h)$ of dimension at least two.
The trace construction gives the Green representation
$\phi=\lambda\int_DG_D(\,\cdot\,,z)\phi(z)\,M_h(dz)$
for each eigenfunction. The bound on $K_2$ makes this representative
bounded. The Frostman estimate also controls the near-diagonal
$G_D^2$ integral uniformly, so a near/far split makes
$x\mapsto G_D(x,\cdot)$ continuous in $L^2(M_h)$.
Thus the eigenfunctions are continuous, as in
\cite[Lemma 6.27]{Cai}. For an orthonormal pair
$\phi,\psi\in V_\lambda$ and a bounded smooth direction $g$, define
\begin{equation}\label{eq:two-response-definition}
 r_g(\phi,\psi)=
 \begin{pmatrix}
  \dfrac12\displaystyle\int_Dg(\phi^2-\psi^2)\,dM_h\\[2mm]
  \displaystyle\int_Dg\phi\psi\,dM_h
 \end{pmatrix}\in\R^2.
\end{equation}

\begin{lemma}\label{lem:ranktwo}
For every orthonormal $\phi,\psi\in V_\lambda$,
\begin{equation}\label{eq:rankmap}
 \operatorname{span}\{r_g(\phi,\psi):
                         g\in C_c^\infty(D)\}=\R^2.
\end{equation}
There exist $M\ge2$ and $c_\lambda>0$ such that
\begin{equation}\label{eq:ranktwo-uniform}
 \det\left(\sum_{l=1}^M
       r_{e_l}(\phi,\psi)r_{e_l}(\phi,\psi)^T\right)
 \ge c_\lambda
\end{equation}
for every orthonormal pair in $V_\lambda$.
\end{lemma}
\begin{proof}
If \eqref{eq:rankmap} fails, some $(a,b)\ne(0,0)$ satisfies
\[
 \int_Dg\bigl[a(\phi^2-\psi^2)+2b\phi\psi\bigr]\,dM_h=0
 \qquad\text{for every }g\in C_c^\infty(D).
\]
The continuous function in brackets vanishes everywhere, because
$M_h$ has full support. An orthogonal diagonalization of
$\left(\begin{smallmatrix}a&b\\b&-a\end{smallmatrix}\right)$,
whose eigenvalues are $\pm\sqrt{a^2+b^2}$, gives orthonormal
$u,v\in V_\lambda$ with
\begin{equation}\label{eq:rank-contradiction}
 u^2-v^2=(u+v)(u-v)=0\quad\text{on }D.
\end{equation}
Since $u+v$ is not identically zero, it is nonzero on some ball
$B\Subset D$. Thus $u-v=0$ on $B$. Lemma~\ref{lem:ucp} implies
$u-v=0$ on $D$, contradicting orthonormality.

We next pass to the fixed basis $(e_l)$. Boundedness of the
eigenfunctions and \eqref{eq:green-frostman-potential} imply finite
Green energy of $(\phi^2-\psi^2)M_h$ and $\phi\psi M_h$.
For either density $\theta$, Green duality gives
\begin{equation}\label{eq:response-direction-continuity}
 \left|\int_Dg\theta\,dM_h\right|
 \le\|g\|_H I_G(\theta M_h)^{1/2}.
\end{equation}
Choose $g_1,g_2\in C_c^\infty(D)$ with
$\det(r_{g_1}(\phi,\psi),r_{g_2}(\phi,\psi))\ne0$.
Approximating these directions by their finite expansions in
$(e_l)$ preserves this nonzero determinant by
\eqref{eq:response-direction-continuity}. Some finite collection
$r_{e_1}(\phi,\psi),\ldots,r_{e_M}(\phi,\psi)$ therefore spans
$\R^2$.

The orthonormal two-frame space of $V_\lambda$ is compact.
The determinants in \eqref{eq:ranktwo-uniform} are continuous in
the frame and increase with $M$. Indeed, writing
$r_l=r_{e_l}(\phi,\psi)$, the Cauchy--Binet identity gives
\[
 \det\left(\sum_{l=1}^Mr_lr_l^T\right)
   =\sum_{1\le l<m\le M}\det(r_l,r_m)^2.
\]
Their positivity sets form an increasing open cover of the frame
space. A finite subcover gives a common $M$, and the minimum at
that $M$ is strictly positive. This proves
\eqref{eq:ranktwo-uniform}.
\end{proof}

\subsubsection{Countable spectral windows and null projections}
Fix $k$. For each $M\ge\max\{2,k\}$, disintegrate the field as
\begin{equation}\label{eq:fiber-disintegration}
 h=\eta+\sum_{l=1}^Mt_le_l,\qquad
 \PP(dh)=\PP_M^\perp(d\eta)\,\varphi_M(t)\,dt,\qquad
 \varphi_M(t)=(2\pi)^{-M/2}e^{-|t|^2/2}.
\end{equation}
For a compact positive operator $A$ and rational $0<a<b$ outside
its spectrum, write
\[
 P_{a,b}(A)=\one_{(a,b)}(A),\qquad
 V_{a,b}(A)=\operatorname{Ran}P_{a,b}(A).
\]
This range is finite dimensional. For an $H$-orthonormal pair
$u,v\in V_{a,b}(T(h))$, define
\begin{equation}\label{eq:two-response-operator}
 r_l(h;u,v)=
 \begin{pmatrix}
 \frac12\bigl(\langle u,\partial_lT(h)u\rangle_H
             -\langle v,\partial_lT(h)v\rangle_H\bigr)\\[1mm]
 \langle u,\partial_lT(h)v\rangle_H
 \end{pmatrix},
 \qquad 1\le l\le M.
\end{equation}
On windows of rank at least two, let
\begin{equation}\label{eq:frame-determinant}
 d_M(h;a,b)=
 \min_{\substack{u,v\in V_{a,b}(T(h))\\
      \|u\|_H=\|v\|_H=1,\ \langle u,v\rangle_H=0}}
 \det\left(\sum_{l=1}^Mr_l(h;u,v)r_l(h;u,v)^T\right).
\end{equation}
The minimum exists by compactness of the two-frame space.

\begin{lemma}\label{lem:generic-lines}
For every fixed $k$, $\PP^{\perp k}$-almost every $h^{\perp k}$
satisfies
\begin{equation}\label{eq:generic-lines-conclusion}
 \lambda_i(h^{\perp k}+te_k)
   <\lambda_{i+1}(h^{\perp k}+te_k)
 \qquad\text{for all }t\in\R,\ i\ge1.
\end{equation}
Consequently the functions
\begin{equation}\label{eq:generic-lines-analytic}
 t\longmapsto Y_i(h^{\perp k}+te_k),\qquad
 t\longmapsto S_n(h^{\perp k}+te_k)
\end{equation}
are locally real analytic on $\R$ for every $i,n\ge1$.
\end{lemma}
\begin{proof}
Let $\Omega_*$ be the $V$-invariant set from
Lemma~\ref{lem:coherent-shifts}. For fixed $M$ put
\[
 h_\eta(t)=\eta+\sum_{l=1}^Mt_le_l,\qquad
 \Omega_*(\eta)=\{t\in\R^M:h_\eta(t)\in\Omega_*\}.
\]
Invariance makes each section either empty or all of $\R^M$.
Since $\PP(\Omega_*)=1$, Gaussian disintegration gives
\begin{equation}\label{eq:full-fibers}
 1=\int\left(\int_{\R^M}\one_{\Omega_*(\eta)}(t)
                          \varphi_M(t)\,dt\right)\PP_M^\perp(d\eta),
 \qquad
 \Omega_*(\eta)=\R^M\quad\text{for }\PP_M^\perp\text{-a.e. }\eta.
\end{equation}
Fix such an $\eta$. By Lemma~\ref{lem:coherent},
$T_\eta(t):=T(h_\eta(t))$ is compact, positive and real analytic
in Hilbert--Schmidt norm on $\R^M$.

With $M,\eta$ fixed, for rational $0<a<b$ and $j\ge1$ define
\[
 U_{a,b,j}=
 \left\{t\in\R^M:
 \begin{array}{l}
 a,b\notin\operatorname{spec}T_\eta(t),\\
 \dim V_{a,b}(T_\eta(t))\ge2,\\
 d_M(h_\eta(t);a,b)>1/j
 \end{array}\right\}.
\]
Near a point where the endpoints stay outside the spectrum, the
Riesz projection has fixed rank and a continuous orthonormal frame.
The minimum in \eqref{eq:frame-determinant} is then taken over a
fixed compact frame space and varies continuously. Thus
$U_{a,b,j}$ is open.

These windows detect every repeated positive eigenvalue when
$M,a,b,j$ vary. Indeed, if $\tau>0$ is a repeated eigenvalue of
$T(h)$, choose rational endpoints such that
\begin{equation}\label{eq:isolating-window}
 \operatorname{spec}T(h)\cap[a,b]=\{\tau\}.
\end{equation}
Set $\lambda=\tau^{-1}$. Every $H$-orthonormal pair in this
eigenspace has the form $u=\phi/\sqrt\lambda$,
$v=\psi/\sqrt\lambda$, with $\phi,\psi$ orthonormal in
$L^2(M_h)$. Differentiating the trace pairing gives
\begin{equation}\label{eq:two-response-identification}
 \begin{aligned}
 \langle u,\partial_lT(h)v\rangle_H
   &=\gamma\int_De_luv\,dM_h
     =\gamma\tau\int_De_l\phi\psi\,dM_h,\\
 r_l(h;u,v)&=\gamma\tau\,r_{e_l}(\phi,\psi).
 \end{aligned}
\end{equation}
Lemma~\ref{lem:ranktwo} supplies an $M\ge\max\{2,k\}$ for which
\begin{equation}\label{eq:window-positive}
 d_M(h;a,b)\ge(\gamma\tau)^4c_\lambda>0.
\end{equation}
Choosing $j$ large enough places the corresponding parameter in
$U_{a,b,j}$. The choices of $(M,a,b,j)$ are countable.

We now bound the dimension of the collision set in one window.
Choose a sufficiently small ball $U\Subset U_{a,b,j}$.
An analytic orthonormal frame of the Riesz range represents the
restriction of $T_\eta(t)$ by a real-analytic symmetric matrix
$H_U(t)\in\operatorname{Sym}(d)$. Its collision set is
\begin{equation}\label{eq:analytic-discriminant}
 Z_U=\{t\in U:\Delta(t)=0\},\qquad
 \Delta(t)=\operatorname{disc}_z\det(zI-H_U(t)).
\end{equation}
The subsets of $\operatorname{Sym}(d)$ with a fixed eigenvalue
multiplicity pattern are semialgebraic, as follows from the rank
conditions for the characteristic polynomial and its derivatives.
Their inverse images under $H_U$ are semianalytic. By
\cite[Proposition 2.10 and Corollary 2.11]{BM}, $Z_U$ therefore has
a locally finite stratification into smooth analytic manifolds on
which these multiplicities are constant.

On one such stratum $S$, fix a repeated branch $\tau(t)$ of
multiplicity $d_0\ge2$. Near $t_0\in S$, choose $r>0$ so that
the positively oriented circle $|z-\tau(t_0)|=r$ encloses only
$\tau(t_0)$ from the spectrum at $t_0$. The analytic projection is
\[
 P_\tau(t)=\frac1{2\pi i}\oint_{|z-\tau(t_0)|=r}
                  (z-T_\eta(t))^{-1}\,dz.
\]
It has rank $d_0$, and on $S$,
$\tau(t)=d_0^{-1}\Tr(T_\eta(t)P_\tau(t))$. In particular,
$\tau|_S$ and $P_\tau|_S$ are $C^1$.
For $w\in T_tS$, differentiate $T_\eta P_\tau=\tau P_\tau$
along a curve in $S$ with velocity $w$. Multiplication by
$P_\tau$ on both sides cancels the projection derivatives and gives
\begin{equation}\label{eq:stratum-tangent}
 P_\tau(t)\partial_wT_\eta(t)P_\tau(t)
       =(\partial_w\tau(t))P_\tau(t).
\end{equation}
For an $H$-orthonormal pair $u,v\in\operatorname{Ran}P_\tau(t)$,
define the $2\times M$ matrix
\begin{equation}\label{eq:linear-splitting-map}
 B=\bigl(r_1(h_\eta(t);u,v),\ldots,r_M(h_\eta(t);u,v)\bigr).
\end{equation}
Equation \eqref{eq:stratum-tangent} implies $T_tS\subset\ker B$.
On the other hand, \eqref{eq:frame-determinant} and
$t\in U_{a,b,j}$ give
\begin{equation}\label{eq:rank-codimension}
 \det(BB^T)
 =\det\left(\sum_{l=1}^M
   r_l(h_\eta(t);u,v)r_l(h_\eta(t);u,v)^T\right)>1/j.
\end{equation}
Thus $B$ has rank two, and
\begin{equation}\label{eq:codimension-two}
 \dim S\le\dim\ker B=M-2.
\end{equation}

Let $\pi_k:\R^M\to\R^{M-1}$ delete the $k$th coordinate.
Each stratum has a countable atlas of locally Lipschitz charts
of dimension at most $M-2$. Hence
\begin{equation}\label{eq:null-projection}
 \mathcal H^{M-1}(\pi_kS)=0,\qquad
 \operatorname{Leb}_{M-1}(\pi_kZ_U)=0.
\end{equation}
Countably many balls $U$ with the above Riesz representation cover
each $U_{a,b,j}$. Define
\begin{equation}\label{eq:good-collision-set}
 Z_{M,\eta}=
 \bigcup_{\substack{a,b\in\mathbb Q,\ 0<a<b\\j\ge1}}
 \left\{t\in U_{a,b,j}:
     T_\eta(t)\text{ has a repeated eigenvalue in }(a,b)\right\}.
\end{equation}
Countability and \eqref{eq:null-projection} imply
\begin{equation}\label{eq:projection-collision}
 \operatorname{Leb}_{M-1}(\pi_kZ_{M,\eta})=0.
\end{equation}

To integrate this statement in $\eta$, realize the complementary
field by its remaining Gaussian coordinates. On $\Omega_*$,
the finite matrix truncations of $T(h)$ converge in
$\mathrm{HS}(H)$, and their coefficients are Borel GMC integrals.
Extending $T$ by zero off $\Omega_*$ therefore gives a Borel
operator map. Each $\partial_lT$ is a Hilbert--Schmidt
finite-difference limit on $\Omega_*$ and is Borel there as well.
The spectral projections and the compact-frame minima in
\eqref{eq:frame-determinant} are Borel on separated windows.
Thus the joint set $\{(\eta,t):t\in Z_{M,\eta}\}$ is Borel,
taking empty sections when $\eta\notin\Omega_*$.
Its projection under $(\eta,t)\mapsto(\eta,t_{\widehat k})$ is
analytic and therefore universally measurable. Gaussian Fubini
applies on the completed product space, and
\eqref{eq:projection-collision} gives
\begin{equation}\label{eq:line-null-disintegration}
 \int\!\int_{\R^{M-1}}
  \one_{\pi_kZ_{M,\eta}}(t_{\widehat k})
  \varphi_{M-1}(t_{\widehat k})\,dt_{\widehat k}\,
  \PP_M^\perp(d\eta)=0.
\end{equation}
Under this product measure,
$\eta+\sum_{l\le M,\ l\ne k}t_le_l$ has law $\PP^{\perp k}$.
Consequently, for each $M$, almost every $e_k$-line misses the
collision set detected at that $M$. Intersect these full-measure
sets over $M\ge\max\{2,k\}$. Every collision on a coherent line
is detected for some $M$ by \eqref{eq:isolating-window}--%
\eqref{eq:window-positive}; hence no collision remains on any
line in this intersection. This proves
\eqref{eq:generic-lines-conclusion} and \eqref{eq:linewise-target}.
Each eigenvalue is simple and isolated at every parameter on these
lines, so analytic perturbation yields
\eqref{eq:generic-lines-analytic}.
\end{proof}

\subsubsection{Second response and coordinate integration by parts}
Fix $g\in V$. At a simple spectrum, choose the locally analytic
branches along $h+tg$, normalized in $L^2(M_{h+tg})$, and put
$m_{im}(g)=\int_Dg\phi_i\phi_m\,dM_h$.
Write $c_m=\int_D\dot\phi_i\phi_m\,dM_h$, with a dot denoting
the derivative at $t=0$. Differentiating the weak eigenvalue
equation against $\phi_m$ and differentiating the normalization give
\[
 (\lambda_m-\lambda_i)c_m=\gamma\lambda_i m_{im}(g)
       \quad(m\ne i),\qquad
 c_i=-\frac\gamma2m_{ii}(g).
\]
The first response \eqref{eq:firstresponse} is
$\partial_gY_i=-\gamma m_{ii}(g)$. To differentiate it once more,
use Parseval for $g\phi_i$:
\[
 \begin{aligned}
 \partial_gm_{ii}(g)
 &=2\sum_m c_m m_{im}(g)
                  +\gamma\int_Dg^2\phi_i^2\,dM_h\\
 &=\gamma\sum_{m\ne i}
       \left(1+\frac{2\lambda_i}{\lambda_m-\lambda_i}\right)
                                                 m_{im}(g)^2.
 \end{aligned}
\]
Here $2c_i m_{ii}(g)=-\gamma m_{ii}(g)^2$ cancels the diagonal
term from $\gamma\int_Dg^2\phi_i^2\,dM_h$.
It follows that
\begin{equation}\label{eq:secondresponse}
 \partial_g^2Y_i
 =\gamma^2\sum_{m\ne i}
       \frac{\lambda_i+\lambda_m}{\lambda_i-\lambda_m}m_{im}(g)^2.
\end{equation}
For fixed $g$ the series is absolutely convergent: the ratio is
bounded for large $m$, the remaining denominators are nonzero,
and $\sum_m m_{im}(g)^2=\int_Dg^2\phi_i^2\,dM_h<\infty$.

When \eqref{eq:secondresponse} is summed over $i\le n$, the terms
with $i,m\le n$ cancel in pairs because $m_{im}(g)=m_{mi}(g)$.
Only pairs crossing the $n$th spectral level remain:
\begin{equation}\label{eq:secondcumulative}
 \partial_g^2S_n
 =-\gamma^2\sum_{i\le n<m}
       \frac{\lambda_m+\lambda_i}{\lambda_m-\lambda_i}m_{im}(g)^2
 \le0.
\end{equation}
The product min--max principle extends concavity across collisions:
\[
 S_n(h+tg)=\inf_{\substack{L\subset H\\\dim L=n}}
 \left\{\log\det(\langle v_a,v_b\rangle_H)
       -\log\det\left(\int_Dv_av_b e^{\gamma tg}\,dM_h\right)\right\},
\]
where $(v_a)_{a=1}^n$ is any basis of $L$ and $-\log0=+\infty$.
The mass Gram determinant equals
\[
 \frac1{n!}\int_{D^n}\det(v_a(x_b))^2
       e^{\gamma t\sum_b g(x_b)}\prod_{b=1}^nM_h(dx_b).
\]
Its logarithm is convex by H\"older's inequality. Each expression
in the infimum is therefore concave in $t$, and so is $S_n(h+tg)$.

\begin{corollary}\label{cor:coordinate-ibp}
For each $n,k$, the Gaussian-coordinate distributional second
derivative of $S_n$ has no singular part and belongs to $L^1(\PP)$.
For $F=u(Y)$ with $u\in C_c^\infty(\R^N)$, and for Gaussian
cylinder tests with a $C_b^\infty$ profile,
\begin{equation}\label{eq:coordinate-ibp}
 \EE[\partial_kS_n\,\partial_kF]
 =-\EE[(\partial_k^2S_n-\xi_k\partial_kS_n)F].
\end{equation}
Both sides are integrable for each fixed $k$.
\end{corollary}
\begin{proof}
On almost every $e_k$-line, Lemma~\ref{lem:generic-lines} gives
analyticity on $\R$. By \eqref{eq:secondcumulative} and
\eqref{eq:direction-lipschitz}, the function $S_n$ on this line
satisfies
\[
 S_n''\le0,\qquad |S_n'|\le L_k:=n\gamma\|e_k\|_\infty.
\]
Let $\varphi_1(t)=(2\pi)^{-1/2}e^{-t^2/2}$. Choose smooth compactly
supported cutoffs $\chi_R\uparrow1$ with
$\|\chi_R'\|_\infty\le C/R$. Integration by parts gives
\[
 \int_\R\chi_R\varphi_1(t)(-S_n'')\,dt
 =-\int_\R t\chi_R\varphi_1(t) S_n'\,dt
       +\int_\R\chi_R'\varphi_1(t) S_n'\,dt.
\]
The last term tends to zero. Monotone convergence on the left and
dominated convergence on the right yield
\[
 \int_\R\varphi_1(t)(-S_n'')\,dt
 =-\int_\R t\varphi_1(t) S_n'\,dt
 \le L_k\int_\R|t|\varphi_1(t)\,dt.
\]
The bound is uniform in the complementary field, so
$\EE|\partial_k^2S_n|\le L_k\EE|\xi_k|$.
Also $|\xi_k\partial_kS_n|\le L_k|\xi_k|$.
For a spectral test, boundedness of $\nabla u$ and
\eqref{eq:direction-lipschitz} make $F$ bounded and Lipschitz on
each line; the same holds for the stated cylinder tests.
Weighted integration by parts on the line, followed by Gaussian
disintegration, now gives \eqref{eq:coordinate-ibp}.
The line derivatives are classical everywhere on these
full-measure fibers, which also identifies the distributional
second derivative.
\end{proof}

\subsection{Operator-valued Gaussian calculus}
\label{app:operator}
The coordinate identities of Appendix~\ref{app:curvature} will be
summed through the Gaussian generator of $T$. We first construct
its Malliavin derivative tensors and estimate their moments, then
pass to the generator and its finite-coordinate approximations.
These estimates enter the compensation identity in
Appendix~\ref{app:drift}.

Set $m_h=M_h(D)$ and $\bar M_h=m_h^{-1}M_h$.
Since $C_\alpha\ge m_h$, the ratio $A=C_\alpha/m_h$ is at least
one, and \eqref{eq:Frostman} gives
$\bar M_h(B(x,r))\le\min\{1,Ar^\alpha\}$ for $0<r\le1$.
For fixed $x$, the tail of $\log^+(1/|x-z|)$ under $\bar M_h$
is therefore bounded by $\min\{1,Ae^{-\alpha t}\}$.
Splitting the layer-cake integral at $t=\alpha^{-1}\log A$
and using the logarithmic Green bound yields, for every integer
$l\ge1$,
\begin{equation}\label{eq:masslog}
 \iint_{D^2}G_D(x,z)^l\,M_h(dx)M_h(dz)
 \le C_lm_h^2\left(1+\log\frac{C_\alpha}{m_h}\right)^l.
\end{equation}
Here and below $C_l$ may also depend on the fixed $D,\alpha$.

The logarithmic factor has moments of every order.
Indeed, choose $0<\theta<q$ and conjugate exponents $a,a'>1$
such that $\theta a<q$. H\"older's inequality gives
\begin{equation}\label{eq:massratio-moment}
 \EE\left(\frac{C_\alpha}{m_h}\right)^\theta
 \le(\EE C_\alpha^{\theta a})^{1/a}
       (\EE m_h^{-\theta a'})^{1/a'}<\infty.
\end{equation}
The first factor is finite by Lemma~\ref{lem:frostman}.
For the second, $m_h\ge M_h(B)$ on any interior ball
$B\Subset D$, whose mass has all negative moments
\cite[Theorem 2.12]{RV}. Since
$(1+\log A)^v\le C_{\theta,v}A^\theta$ for $A\ge1$,
\begin{equation}\label{eq:massratio-logmoments}
 \EE\left(1+\log\frac{C_\alpha}{m_h}\right)^v<\infty,
 \qquad 0<v<\infty.
\end{equation}

For $g_1,\ldots,g_r\in V$ define the response operators by
\begin{equation}\label{eq:T-diff}
 \langle T_r[g_1,\ldots,g_r]u,v\rangle_H
 =\gamma^r\int_Dg_1\cdots g_r uv\,dM_h,
 \qquad u,v\in H,\quad r\ge0,
\end{equation}
with $T_0=T$. These are Hilbert--Schmidt operators by the
bounded-multiplier estimate in Lemma~\ref{lem:coherent}.
We now extend the directional prescription to Hilbert tensors
and identify their closed Gaussian derivatives.

\begin{lemma}\label{lem:T-derivatives}
For every integer $r\ge0$ and $1<p<4/\gamma^2$,
\eqref{eq:T-diff} extends to a tensor
$T_r\in L^p(\PP;H^{\otimes r}\otimes\HS(H))$.
Moreover, $T\in\mathbb D^{r,p}(\HS(H))$, $\D^rT=T_r$, and
\begin{equation}\label{eq:T-Sobolev}
 \|T_r\|^2
 =\gamma^{2r}\iint_{D^2}G_D(x,z)^{r+2}\,
                              M_h(dx)M_h(dz).
\end{equation}
\end{lemma}
\begin{proof}
Fix $r$. Under the isometric identification
$\HS(H)\simeq H\otimes H$, the proposed tensor has coefficients
\[
 c_{\mathbf k}=\gamma^r\int_D
              \prod_{j=1}^{r+2}e_{k_j}(x)\,M_h(dx),
 \qquad \mathbf k=(k_1,\ldots,k_{r+2}).
\]
They are finite because each $e_k$ is bounded. To estimate their
square sum, damp all $r+2$ indices by the heat semigroup:
\[
 c_{\mathbf k}^{\,\varepsilon}
 =c_{\mathbf k}\prod_{j=1}^{r+2}e^{-\varepsilon\kappa_{k_j}}.
\]
The resulting Green kernel is
\begin{equation}\label{eq:regularized-green-tensor}
 G_\varepsilon(x,z)
 =\sum_{k\ge1}e^{-2\varepsilon\kappa_k}e_k(x)e_k(z)
 =2\pi\int_{2\varepsilon}^\infty p_t^D(x,z)\,dt,
\end{equation}
where $p_t^D$ is the kernel of $P_t^D=e^{t\Delta_D}$.
Thus $G_\varepsilon\ge0$ and
$G_\varepsilon\uparrow G_D$ off the diagonal as
$\varepsilon\downarrow0$. For fixed $\varepsilon>0$,
Cauchy--Schwarz and the heat-kernel bound give
\[
 \sum_{k\ge1}e^{-2\varepsilon\kappa_k}|e_k(x)e_k(z)|
 \le\bigl(G_\varepsilon(x,x)G_\varepsilon(z,z)\bigr)^{1/2}
 \le C_\varepsilon.
\]
This absolute bound and finiteness of $M_h$ justify expanding
the product of $r+2$ kernels and exchanging sums with integrals.
Consequently
\begin{equation}\label{eq:heat-tensor-parseval}
 \begin{aligned}
 \sum_{\mathbf k}|c_{\mathbf k}^{\,\varepsilon}|^2
 &=\gamma^{2r}\iint_{D^2}G_\varepsilon(x,z)^{r+2}\,
                                  M_h(dx)M_h(dz)\\
 &\le\gamma^{2r}\iint_{D^2}G_D(x,z)^{r+2}\,
                                  M_h(dx)M_h(dz)<\infty.
 \end{aligned}
\end{equation}
The last bound is \eqref{eq:masslog}.
The squared coefficient weights increase to one, while the
nonnegative kernels increase to $G_D$. The diagonal is
$M_h\otimes M_h$-null by \eqref{eq:Frostman}.
Monotone convergence on both sides of
\eqref{eq:heat-tensor-parseval} therefore gives
\begin{equation}\label{eq:tensor-norm-limit}
 \sum_{\mathbf k}|c_{\mathbf k}|^2
 =\gamma^{2r}\iint_{D^2}G_D(x,z)^{r+2}\,
                                  M_h(dx)M_h(dz)<\infty.
\end{equation}
This constructs $T_r\in H^{\otimes r}\otimes\HS(H)$ with
the norm in \eqref{eq:T-Sobolev}. Its coefficients agree with
\eqref{eq:T-diff} for directions in $V$; density then extends
the response continuously to all directions in $H$.

We next prove the $L^p$ bound. Choose $p<p_1<4/\gamma^2$.
A finite cover $\mathcal Q_0$ of $D$ by unit dyadic squares and
\eqref{eq:dyadic-mass-moment} give
\[
 \EE m_h^{p_1}
 \le(\#\mathcal Q_0)^{p_1-1}
       \sum_{Q\in\mathcal Q_0}\EE M_h(Q\cap D)^{p_1}<\infty.
\]
Combining \eqref{eq:masslog} with H\"older's inequality and
\eqref{eq:massratio-logmoments}, we obtain
\begin{equation}\label{eq:Malliavin-tensor-moments}
 \begin{aligned}
 \EE\|T_r\|^p
 &\le C_{r,p}\EE\left[
       m_h^p\left(1+\log\frac{C_\alpha}{m_h}\right)^{(r+2)p/2}
                         \right]\\
 &\le C_{r,p}(\EE m_h^{p_1})^{p/p_1}
       \left[\EE\left(1+\log\frac{C_\alpha}{m_h}\right)^{
                    \frac{(r+2)pp_1}{2(p_1-p)}}\right]^{1-p/p_1}
 <\infty.
 \end{aligned}
\end{equation}

It remains to identify these tensors as Malliavin derivatives.
For $g_1,\ldots,g_j,u,v\in V$, the coherent-shift identity
\eqref{eq:alltilts} gives
\begin{equation}\label{eq:tensor-directional-derivative}
 \left.\frac{d}{dt}\right|_{t=0}
   \langle T_j(h+te_k)[g_1,\ldots,g_j]u,v\rangle_H
 =\langle T_{j+1}[e_k,g_1,\ldots,g_j]u,v\rangle_H.
\end{equation}
Let $\mathcal G_m=\sigma(\xi_1,\ldots,\xi_m)$ and
$T^{(m)}=\EE[T\mid\mathcal G_m]$.
For $0\le j<r$ and a deterministic finite tensor
$A\in H^{\otimes j}\otimes\HS(H)$ with indices in the fixed basis,
ordinary Gaussian integration by parts on the $k$th coordinate
line yields
\[
 \EE[\langle T_j,A\rangle\,\partial_k\Phi]
 =\EE\bigl[
   \bigl(\xi_k\langle T_j,A\rangle
       -\langle T_{j+1},e_k\otimes A\rangle\bigr)\Phi\bigr],
\]
where $\Phi=\varphi(\xi_1,\ldots,\xi_m)$,
$\varphi\in C_c^\infty(\R^m)$, and $k\le m$.
The moment bound \eqref{eq:Malliavin-tensor-moments} justifies
disintegration and passage from finite tensors to arbitrary tests
in the tensor Hilbert space. Let $\Pi_m$ project onto
$\operatorname{span}\{e_1,\ldots,e_m\}$.
Conditioning and iterating give weak derivatives of $T^{(m)}$
through order $r$, with the same $L^p$ bounds.
Finite projections of the Hilbert--Schmidt range, followed by
compact cutoffs and mollification in the $m$ Gaussian variables,
therefore give $T^{(m)}\in\mathbb D^{r,p}(\HS(H))$, with
\begin{equation}\label{eq:finite-Gaussian-tensor}
 \D^jT^{(m)}
 =\Pi_m^{\otimes j}\EE[T_j\mid\mathcal G_m],
 \qquad 0\le j\le r,
\end{equation}
The contraction properties of $\Pi_m^{\otimes j}$ and conditional
expectation give
\begin{equation}\label{eq:conditional-tensor-convergence}
 \begin{aligned}
 \|\D^jT^{(m)}-T_j\|_{L^p}
 &\le\|\EE[T_j\mid\mathcal G_m]-T_j\|_{L^p}\\
 &\quad+\|(\Pi_m^{\otimes j}-I)T_j\|_{L^p}
 \longrightarrow0.
 \end{aligned}
\end{equation}
The first term tends to zero by martingale convergence, and the
second by strong convergence of the projections and dominated
convergence. For $j=0$ the same argument gives $T^{(m)}\to T$.
Closedness of the Malliavin derivatives now proves
$T\in\mathbb D^{r,p}(\HS(H))$ and $\D^jT=T_j$ for $j\le r$.
\end{proof}

\begin{proposition}\label{prop:T-Sobolev}
For $1<p<4/\gamma^2$ and every $s\ge0$,
\begin{equation}\label{eq:operator-generator-definition}
 \LL_sT=-\delta(Q_s\D T)\in L^p(\PP;\HS(H)).
\end{equation}
The finite-coordinate generators are
\begin{equation}\label{eq:operator-coordinate-generator}
 \LL_{s,K}T:=-\delta(\Pi_KQ_s\D T)
 =\sum_{k=1}^Kq_k(s)(\partial_k^2T-\xi_k\partial_kT).
\end{equation}
They satisfy
\begin{align}
 \LL_{s,K}T&\longrightarrow\LL_sT
     &&\text{in }L^p(\PP;\HS(H)),\qquad K\to\infty,
   \label{eq:LKconv}\\
 \LL_sT&\longrightarrow\LL_0T
     &&\text{in }L^p(\PP;\HS(H)),\qquad s\downarrow0,
   \label{eq:Lsconv}
\end{align}
and
\begin{equation}\label{eq:OU-uniform-bound}
 \sup_{K\ge1,\ s\ge0}
        \|\LL_{s,K}T\|_{L^p(\PP;\HS(H))}<\infty.
\end{equation}
\end{proposition}
\begin{proof}
Lemma~\ref{lem:T-derivatives} gives $T_1=\D T$ in
$\mathbb D^{1,p}(H\otimes\HS(H))$, with derivative $T_2$.
For a deterministic bounded operator $B:H\to H$,
\begin{equation}\label{eq:weighted-operator-derivative}
 \D(BT_1)=(I_H\otimes B)T_2,
\end{equation}
where the new derivative slot comes first and the identity on
$\HS(H)$ is implicit.
The Hilbert-valued divergence estimate
\cite[Proposition 4.3]{PV} therefore gives
\begin{equation}\label{eq:Gaussian-divergence-bound}
 \|\delta(BT_1)\|_{L^p(\PP;\HS(H))}
 \le C_p\bigl(\|BT_1\|_{L^p}
               +\|(I_H\otimes B)T_2\|_{L^p}\bigr).
\end{equation}
Taking $B=Q_s$ and $B=\Pi_KQ_s$, both contractions, proves
\eqref{eq:operator-generator-definition} and
\eqref{eq:OU-uniform-bound}.
Taking $B=(\Pi_K-I)Q_s$ or $B=Q_s-I$ proves the two limits:
these operators converge strongly to zero, are uniformly bounded,
and hence send the two terms on the right to zero by dominated
convergence in the Hilbert tensor spaces.

For the coordinate formula, $\xi_k\partial_kT\in L^p$ by
H\"older's inequality and Lemma~\ref{lem:T-derivatives} at any
$p_1\in(p,4/\gamma^2)$. The Gaussian product rule
$\delta(Fe_k)=\xi_kF-\partial_kF$, applied to
$F=\partial_kT$, gives \eqref{eq:operator-coordinate-generator}.
\end{proof}

\subsection{Compensated spectral drift}
\label{app:drift}
We combine the coordinate integration-by-parts identity of
Appendix~\ref{app:curvature} with the operator estimates of
Appendix~\ref{app:operator}. At a finite Gaussian cutoff, the
internal spectral pairs cancel their denominators. The remaining
nonnegative crossing flux is controlled by taking expectations
in the resulting compensation identity.

\subsubsection{Finite-coordinate compensation}
Put $C_n=\sum_{i\le n}\lambda_i(v_i\otimes v_i)$ for $n\ge1$.
For $i,m,k\ge1$, write
\begin{equation}\label{eq:drift-coordinate-products}
 m_{im,k}=m_{im}(e_k)=\int_De_k\phi_i\phi_m\,dM_h,
 \qquad d_{i,k}=\int_De_k^2\,d\nu_i.
\end{equation}
Parseval in $L^2(M_h)$ and $\nu_i(D)=1$ give
\begin{equation}\label{eq:drift-parseval}
 \sum_{m\ge1}m_{im,k}^2=d_{i,k}\le\|e_k\|_\infty^2.
\end{equation}
For $s\ge0$ and $K\ge1$, define
\begin{equation}\label{eq:drift-truncations}
 \begin{aligned}
 I_{im,s}^{[K]}&=\sum_{k=1}^Kq_k(s)m_{im,k}^2,\\
 F_{n,s}^{[K]}&=2\gamma^2
       \sum_{i\le n<m}\frac{\lambda_i}{\lambda_m-\lambda_i}
                                               I_{im,s}^{[K]},\\
 \beta_{n,s}^{[K]}&=\sum_{k=1}^Kq_k(s)
                    (\partial_k^2S_n-\xi_k\partial_kS_n).
 \end{aligned}
\end{equation}
For $i\le n<m$,
$\lambda_i/(\lambda_m-\lambda_i)
\le\lambda_n/(\lambda_{n+1}-\lambda_n)$.
Thus the spectral series at fixed $K$ is finite almost surely:
\begin{equation}\label{eq:flux-fixed-cutoff-bound}
 0\le F_{n,s}^{[K]}
 \le2\gamma^2\frac{\lambda_n}{\lambda_{n+1}-\lambda_n}
                      \sum_{i\le n}\sum_{k\le K}d_{i,k}<\infty.
\end{equation}

\begin{lemma}\label{lem:compensated}
For every $n,K\ge1$ and $s\ge0$,
\begin{equation}\label{eq:finite-drift}
 \beta_{n,s}^{[K]}
 =-\Tr(C_n\LL_{s,K}T)
       +\gamma^2\sum_{i,j\le n}I_{ij,s}^{[K]}
       -F_{n,s}^{[K]}.
\end{equation}
Moreover, $\beta_{n,s}^{[K]}\in L^1(\PP)$ and
$\EE\beta_{n,s}^{[K]}=0$.
\end{lemma}
\begin{proof}
Fix $k$ and abbreviate $m_{im}=m_{im,k}$.
The response formula \eqref{eq:T-diff}, evaluated at
$v_i=\lambda_i^{-1/2}\phi_i$, gives
\begin{equation}\label{eq:drift-operator-diagonal}
 \lambda_i\langle v_i,
          (\partial_k^2-\xi_k\partial_k)Tv_i\rangle_H
       =\gamma^2d_{i,k}-\gamma\xi_km_{ii}.
\end{equation}
In the second response \eqref{eq:secondresponse}, use
\[
 \frac{\lambda_i+\lambda_m}{\lambda_i-\lambda_m}
   =-1+\frac{2\lambda_i}{\lambda_i-\lambda_m},
 \qquad
 \sum_{m\ne i}m_{im}^2=d_{i,k}-m_{ii}^2.
\]
Together with $\partial_kY_i=-\gamma m_{ii}$ and
\eqref{eq:drift-operator-diagonal}, this yields
\begin{equation}\label{eq:drift-coordinate-expansion}
 \begin{aligned}
 (\partial_k^2-\xi_k\partial_k)Y_i
 &=-\lambda_i\langle v_i,
           (\partial_k^2-\xi_k\partial_k)Tv_i\rangle_H
       +\gamma^2m_{ii}^2\\
 &\quad+2\gamma^2\lambda_i
             \sum_{m\ne i}\frac{m_{im}^2}{\lambda_i-\lambda_m}.
 \end{aligned}
\end{equation}
For $i<j\le n$, symmetry of $m_{ij}$ and
\begin{equation}\label{eq:drift-internal-pairs}
 \frac{\lambda_i}{\lambda_i-\lambda_j}
       +\frac{\lambda_j}{\lambda_j-\lambda_i}=1
\end{equation}
combine the internal terms into
$2\gamma^2\sum_{i<j\le n}m_{ij}^2$.
Together with the diagonal term this is
$\gamma^2\sum_{i,j\le n}m_{ij}^2$.
For $i\le n<m$, only the $i$th term occurs, with a negative
denominator. Summing \eqref{eq:drift-coordinate-expansion}
over $i\le n$ therefore gives
\begin{equation}\label{eq:drift-coordinate-compensated}
 \begin{aligned}
 (\partial_k^2-\xi_k\partial_k)S_n
 &=-\Tr\bigl(C_n(\partial_k^2-\xi_k\partial_k)T\bigr)
       +\gamma^2\sum_{i,j\le n}m_{ij}^2\\
 &\quad-2\gamma^2\sum_{i\le n<m}
               \frac{\lambda_i}{\lambda_m-\lambda_i}m_{im}^2.
 \end{aligned}
\end{equation}
Multiply by $q_k(s)$ and sum over $k\le K$.
The generator formula \eqref{eq:operator-coordinate-generator}
then proves \eqref{eq:finite-drift}.

By Corollary~\ref{cor:coordinate-ibp}, each
$\partial_k^2S_n-\xi_k\partial_kS_n$ is integrable.
Taking the constant test function in that corollary gives its
zero mean. The finite sum in \eqref{eq:drift-truncations}
therefore has the asserted integrability and mean.
\end{proof}

\subsubsection{Uniform flux control and removal of cutoffs}
Define the full pair energies and crossing flux by
\begin{align}
 I_{im,s}&=\sum_{k\ge1}q_k(s)m_{im,k}^2,
          \label{eq:Jdef}\\
 F_{n,s}&=2\gamma^2\sum_{i\le n<m}
           \frac{\lambda_i}{\lambda_m-\lambda_i}I_{im,s}.
          \label{eq:crossing-flux-def}
\end{align}
All summands are nonnegative, so these are the increasing limits
of $I_{im,s}^{[K]}$ and $F_{n,s}^{[K]}$ as $K\to\infty$.

To control the internal energies, put
$\eta_{im}=\phi_i\phi_mM_h$.
The bound $|\eta_{im}|\le(\nu_i+\nu_m)/2$ gives finite Green
energy. Green--Riesz duality and Parseval in the basis $(e_k)$
then imply
\begin{equation}\label{eq:Jbound}
 \begin{aligned}
 0\le I_{im,s}^{[K]}
 &\le I_{im,s}\le I_{im,0}=I_G(\eta_{im})\\
 &\le\tfrac14I_G(\nu_i+\nu_m)
 \le\tfrac12\{I_G(\nu_i)+I_G(\nu_m)\}.
 \end{aligned}
\end{equation}
Here $G_D\ge0$ bounds the signed energy by that of the total
variation, and Green-energy Cauchy--Schwarz gives the last
inequality. Proposition~\ref{prop:AG} yields
\begin{equation}\label{eq:drift-J-moments}
 \EE I_{im,0}^r<\infty\quad(1\le r<\infty),\qquad
 I_{im,s}^{[K]}\longrightarrow I_{im,s}
             \quad\text{in }L^1(\PP),\quad K\to\infty.
\end{equation}

\begin{proposition}\label{prop:flux-integrable}
For every $n\ge1$,
\begin{equation}\label{eq:drift-flux-conclusion}
 F_{n,0}\in L^1(\PP).
\end{equation}
For $s\ge0$, the full compensated drifts
\begin{align}
 \beta_{n,s}&=-\Tr(C_n\LL_sT)
       +\gamma^2\sum_{i,j\le n}I_{ij,s}-F_{n,s},
       \qquad \beta_{0,s}=0,
       \label{eq:Bdef}\\
 B_{i,s}&=\beta_{i,s}-\beta_{i-1,s},\qquad
 b_{i,s}(y)=\EE[B_{i,s}\mid Y=y]
       \label{eq:bdef}
\end{align}
satisfy $\beta_{n,s},B_{i,s}\in L^1(\PP)$ and
$b_{i,s}\in L^1(\rho)$. Write $b_i=b_{i,0}$. Then
\begin{equation}\label{eq:B-K-limit}
 \|\beta_{n,s}^{[K]}-\beta_{n,s}\|_{L^1(\PP)}
       \longrightarrow0\qquad(K\to\infty).
\end{equation}
As $s\downarrow0$,
\begin{align}
 \|B_{i,s}-B_{i,0}\|_{L^1(\PP)}&\longrightarrow0,
       \label{eq:B-s-limit}\\
 \|b_{i,s}-b_i\|_{L^1(\rho)}&\longrightarrow0.
       \label{eq:b-s-limit}
\end{align}
\end{proposition}
\begin{proof}
Choose $1<p<4/\gamma^2$ and $p'=p/(p-1)$.
Since $\|C_n\|_{\HS}=(\sum_{i\le n}\lambda_i^2)^{1/2}
\le\sqrt n\,\lambda_n$, H\"older's inequality gives
\begin{equation}\label{eq:tracebound}
 \|\Tr(C_n\LL_{s,K}T)\|_{L^1}
 \le\sqrt n\,\|\lambda_n\|_{L^{p'}}
                  \|\LL_{s,K}T\|_{L^p(\PP;\HS(H))}
 \le c_n,
 \qquad s\ge0,\ K\ge1.
\end{equation}
The constant is uniform by Lemma~\ref{lem:frostman} and
\eqref{eq:OU-uniform-bound}.
The other terms in \eqref{eq:finite-drift}, apart from the flux,
are integrable by Lemma~\ref{lem:compensated} and
\eqref{eq:Jbound}; hence $F_{n,s}^{[K]}$ is integrable.
Taking expectations at $s=0$ and using
$\EE\beta_{n,0}^{[K]}=0$, we obtain
\begin{equation}\label{eq:drift-flux-mean-identity}
 \begin{aligned}
 \EE F_{n,0}^{[K]}
 &=-\EE\Tr(C_n\LL_{0,K}T)
             +\gamma^2\sum_{i,j\le n}\EE I_{ij,0}^{[K]}\\
 &\le c_n+\gamma^2\sum_{i,j\le n}\EE I_{ij,0}
 =:c_n'<\infty.
 \end{aligned}
\end{equation}
Monotone convergence now gives
$\EE F_{n,0}=\lim_K\EE F_{n,0}^{[K]}\le c_n'$,
proving \eqref{eq:drift-flux-conclusion}.

Positivity and $q_k(s)\le1$ imply
\begin{equation}\label{eq:drift-flux-domination}
 0\le F_{n,s}^{[K]}\le F_{n,s}\le F_{n,0}.
\end{equation}
Moreover, $F_{n,s}^{[K]}\uparrow F_{n,s}$ as $K\to\infty$,
and $F_{n,s}\uparrow F_{n,0}$ as $s\downarrow0$.
The integrable bound in \eqref{eq:drift-flux-domination} therefore
gives both limits
\begin{equation}\label{eq:drift-flux-two-limits}
 \begin{aligned}
 \|F_{n,s}^{[K]}-F_{n,s}\|_{L^1}&\longrightarrow0,
                                      &&K\to\infty,\\
 \|F_{n,s}-F_{n,0}\|_{L^1}&\longrightarrow0,
                                      &&s\downarrow0.
 \end{aligned}
\end{equation}
The same monotone argument, dominated by \eqref{eq:Jbound},
gives $I_{ij,s}^{[K]}\to I_{ij,s}$ and $I_{ij,s}\to I_{ij,0}$
in $L^1$ for every internal pair.
The trace term also belongs to $L^1$ at the full cutoff by
H\"older's inequality and Proposition~\ref{prop:T-Sobolev}.
Thus all quantities in \eqref{eq:Bdef}--\eqref{eq:bdef}
are integrable.

We can now compare the three terms in the compensated drift:
\begin{equation}\label{eq:drift-beta-estimate}
 \begin{aligned}
 \|\beta_{n,s}^{[K]}-\beta_{n,s}\|_{L^1}
 &\le\sqrt n\,\|\lambda_n\|_{L^{p'}}
              \|\LL_{s,K}T-\LL_sT\|_{L^p(\PP;\HS(H))}\\
 &\quad+\gamma^2\sum_{i,j\le n}
                 \|I_{ij,s}^{[K]}-I_{ij,s}\|_{L^1}
       +\|F_{n,s}^{[K]}-F_{n,s}\|_{L^1}.
 \end{aligned}
\end{equation}
The operator term tends to zero by \eqref{eq:LKconv};
the internal energies and crossing flux converge as above.
This proves \eqref{eq:B-K-limit}.
For $s\downarrow0$, the same estimate uses
$\LL_sT-\LL_0T$, $I_{ij,s}-I_{ij,0}$ and $F_{n,s}-F_{n,0}$.
By \eqref{eq:Lsconv} and the preceding limits,
\begin{equation}\label{eq:drift-beta-limits}
 \|\beta_{n,s}-\beta_{n,0}\|_{L^1}\longrightarrow0.
\end{equation}
The finite-cutoff zero means also pass to the limit:
\begin{equation}\label{eq:drift-zero-mean-limit}
 \EE\beta_{n,s}=0,\qquad
 \EE B_{i,s}=0,\qquad
 \int_Eb_{i,s}\,d\rho=0.
\end{equation}
Finally, conditional expectation is an $L^1$ contraction, so
\begin{equation}\label{eq:drift-coordinate-and-conditional-limits}
 \begin{aligned}
 \|b_{i,s}-b_i\|_{L^1(\rho)}
 &\le\|B_{i,s}-B_{i,0}\|_{L^1(\PP)}\\
 &\le\|\beta_{i,s}-\beta_{i,0}\|_{L^1}
       +\|\beta_{i-1,s}-\beta_{i-1,0}\|_{L^1}
 \longrightarrow0.
 \end{aligned}
\end{equation}
This proves \eqref{eq:B-s-limit} and \eqref{eq:b-s-limit}.
\end{proof}

\subsubsection{Distributional drift and form integration by parts}
Lemma~\ref{lem:MalliavinY} and
Proposition~\ref{prop:flux-integrable} give integrability of the
coefficients in \eqref{eq:static_generator}.
Consequently $\ell_su\in L^1(\rho)$ for
$u\in C_c^\infty(\R^N)$.

\begin{proposition}\label{prop:IBP}
For $s\ge0$, $i\le N$, and $u\in C_c^\infty(\R^N)$,
\begin{equation}\label{eq:IBP}
 \int_E\sum_{j=1}^Na_{ij,s}\partial_ju\,d\rho
       =-\int_Eb_{i,s}u\,d\rho.
\end{equation}
For $u\in C_c^\infty(\R^N)$ and
$v\in\F_s\cap L^\infty(\rho)$,
\begin{equation}\label{eq:drift-weak-generator}
 \E_s(u,v)=-\int_Ev\ell_su\,d\rho.
\end{equation}
\end{proposition}
\begin{proof}
Let $F=u(Y)$ and $1\le n\le N$. Multiplying
\eqref{eq:coordinate-ibp} by $q_k(s)$ and summing over $k\le K$
gives
\begin{equation}\label{eq:drift-coordinate-ibp-finite}
 \EE\langle\Pi_KQ_s\D S_n,\D F\rangle_H
       =-\EE[\beta_{n,s}^{[K]}F].
\end{equation}
Here $\D F=\sum_j(\partial_ju)(Y)\D Y_j$.
The left side converges by strong convergence of $\Pi_KQ_s$
and the $L^2$ bounds on $\D S_n,\D F$.
Since $F$ is bounded, \eqref{eq:B-K-limit} handles the right side.
We obtain
\begin{equation}\label{eq:drift-coordinate-ibp-full}
 \EE\langle Q_s\D S_n,\D(u(Y))\rangle_H
       =-\EE[\beta_{n,s}u(Y)].
\end{equation}
Subtract the identity for $n=i-1$ from that for $n=i$, taking
the former to be zero when $i=1$. Conditioning on $Y$ yields
\begin{equation}\label{eq:drift-divergence-chain}
 \begin{aligned}
 \int_E\sum_ja_{ij,s}\partial_ju\,d\rho
 &=\EE\langle Q_s\D Y_i,\D(u(Y))\rangle_H\\
 &=-\EE[B_{i,s}u(Y)]
 =-\int_Eb_{i,s}u\,d\rho.
 \end{aligned}
\end{equation}
This proves \eqref{eq:IBP}. With $p_N$ extended by zero as in
Appendix~\ref{sec:appendix}, it is the distributional identity
\begin{equation}\label{eq:drift-distributional-divergence}
 \sum_j\partial_j(p_Na_{ij,s})=p_Nb_{i,s}
       \quad\text{in }\mathcal D'(\R^N).
\end{equation}
For $u,v\in\Ccore$, apply \eqref{eq:IBP} to
$v\partial_i u$ and sum over $i$.
Expanding $\partial_j(v\partial_i u)$ gives
\[
 \int_E\sum_{i,j}a_{ij,s}\partial_i u\,\partial_jv\,d\rho
 =-\int_Ev\left(\sum_{i,j}a_{ij,s}\partial_{ij}u
                       +\sum_ib_{i,s}\partial_i u\right)d\rho,
\]
which is \eqref{eq:drift-weak-generator} on the core.

To extend the test function, let $v\in\F_s\cap L^\infty(\rho)$
and choose $v_m\in\Ccore$ converging to $v$ in the form norm
and $\rho$-almost everywhere:
\begin{equation}\label{eq:drift-energy-approximation}
 \|v_m-v\|_{L^2(\rho)}^2
       +\E_s(v_m-v,v_m-v)\longrightarrow0.
\end{equation}
Set $M=\|v\|_\infty$. Choose a smooth $1$-Lipschitz function
$\theta$ with $\theta(0)=0$, equal to the identity on a
neighborhood of $[-M,M]$, and bounded by $M+1$ in absolute value.
Then $w_m:=\theta(v_m)\in\Ccore$,
$\|w_m\|_\infty\le M+1$, and $w_m\to v$ in $L^2(\rho)$ and
almost everywhere.
The chain rule on the core gives
\begin{equation}\label{eq:drift-truncation-energy-estimate}
 \begin{aligned}
 \|D_s(w_m-v)\|_{L^2(\PP;H)}
 &\le\|D_s(v_m-v)\|_{L^2(\PP;H)}\\
 &\quad+\|(\theta'(v_m(Y))-1)D_sv\|_{L^2(\PP;H)}
 \longrightarrow0.
 \end{aligned}
\end{equation}
The second term tends to zero by dominated convergence, since
$\theta'(v_m(Y))\to\theta'(v(Y))=1$ and $|\theta'|\le1$.
Thus $w_m\to v$ in the form norm.
Using this convergence on the left and the uniform bound on
$w_m$ against $\ell_su\in L^1(\rho)$ on the right, we conclude
\begin{equation}\label{eq:drift-extended-form-identity}
 \E_s(u,v)=\lim_m\E_s(u,w_m)
       =-\lim_m\int_Ew_m\ell_su\,d\rho
       =-\int_Ev\ell_su\,d\rho.
\end{equation}
\end{proof}

\begin{proposition}\label{prop:Acoeff}
For every $1\le r<\infty$,
\begin{equation}\label{eq:Acoeff}
 \|a_s-a\|_{L^r(\rho)}\longrightarrow0,\qquad
 \|b_{i,s}-b_i\|_{L^1(\rho)}\longrightarrow0
 \quad(1\le i\le N),\qquad s\downarrow0.
\end{equation}
The coefficients $a_{ij,s},a_{ij},b_{i,s},b_i$ belong to
$L^1(\rho)$.
\end{proposition}
\begin{proof}
Conditional Jensen's inequality and \eqref{eq:covconvergence} give
\begin{equation}\label{eq:limit-proof-coefficients}
 \|a_{ij,s}-a_{ij}\|_{L^r(\rho)}
 \le\|A_{ij,s}-A_{ij,0}\|_{L^r(\PP)}\longrightarrow0.
\end{equation}
The drift convergence is \eqref{eq:b-s-limit}.
Integrability follows from the same Gaussian covariance moments
and Proposition~\ref{prop:flux-integrable}.
\end{proof}

\subsection{The canonical increasing form}
\label{app:mosco}
For $0<s<s'$, the weights in \eqref{eq:Q} satisfy
\begin{equation}\label{eq:mosco-bridge}
 0<q_k(s')\le q_k(s)\le1,\qquad
 q_k(s)\longrightarrow1\quad(s\downarrow0).
\end{equation}
We first identify the derivatives available on the completed
spectral form domains. These derivatives will also distinguish the
increasing limit from the minimal closure at $s=0$.

\begin{theorem}\label{thm:Alimit}
For any $s_m\downarrow0$, define
\begin{equation}\label{eq:Aendpoint}
 \F_* =\left\{u\in\bigcap_m\F_{s_m}:
             \sup_m\E_{s_m}(u,u)<\infty\right\},
 \qquad
 \E_*(u,u)=\sup_m\E_{s_m}(u,u).
\end{equation}
The domain and energy are independent of the sequence $(s_m)$.
The form $(\E_*,\F_*)$ is densely defined, closed, conservative,
and strongly local. The forms $\E_s$ Mosco-converge to $\E_*$. Moreover,
\begin{equation}\label{eq:Agrad}
 \F_0^{\min}\subseteq\F_*
 \subseteq\{u\in L^2(\rho):u(Y)\in\mathbb D^{1,2}\},
 \qquad
 \E_*(u,u)=\EE\|\D(u(Y))\|_H^2,
 \quad u\in\F_*.
\end{equation}
The two energies agree on $\F_0^{\min}$.
\end{theorem}

For $u\in\Ccore$ and $s\ge0$, the Gaussian representation of the
preform is
\begin{equation}\label{eq:two-theorem-interface}
 \E_s(u,u)=\sum_{k\ge1}q_k(s)
       \|\partial_k(u(Y))\|_{L^2(\PP)}^2.
\end{equation}

\begin{lemma}\label{lem:weighted-gradient-domain}
For $s\ge0$ and $u\in\F_s$, the random variable $u(Y)$ has weak
Gaussian derivatives $d_k(u)\in L^2(\PP)$ in every coordinate $e_k$.
They are independent of $s$ whenever $u$ belongs to two form
domains, and
\begin{equation}\label{eq:weighted-domain}
 \E_s(u,u)=\sum_{k\ge1}q_k(s)\|d_k(u)\|_{L^2(\PP)}^2.
\end{equation}
\end{lemma}
\begin{proof}
Choose $u_m\in\Ccore$ converging to $u$ in the form norm:
\begin{equation}\label{eq:weighted-core-approx}
 \|u_m-u\|_{L^2(\rho)}^2
   +\E_s(u_m-u,u_m-u)\longrightarrow0.
\end{equation}
For every $k$, \eqref{eq:two-theorem-interface} gives
\begin{equation}\label{eq:weighted-coordinate-Cauchy}
 \begin{aligned}
 \|\partial_k(u_m(Y))-\partial_k(u_l(Y))\|_{L^2(\PP)}^2
 &\le q_k(s)^{-1}\E_s(u_m-u_l,u_m-u_l),\\
 \|u_m(Y)-u(Y)\|_{L^2(\PP)}
 &=\|u_m-u\|_{L^2(\rho)}.
 \end{aligned}
\end{equation}
Here $q_k(s)>0$ for every fixed $k$ and finite $s$. Since the
Gaussian derivative $\partial_k$ is closed, there exists a unique
$d_k(u)\in L^2(\PP)$ such that
\begin{equation}\label{eq:weighted-coordinate-limit}
 \partial_k(u_m(Y))\longrightarrow d_k(u)
 \quad\text{in }L^2(\PP),\qquad
 d_k(u)=\partial_k(u(Y)).
\end{equation}
The last derivative is understood weakly. In particular, it does
not depend on which $s$ was used to construct the approximation.

The closed gradient $D_s$ constructed in
Proposition~\ref{prop:forms} satisfies
\begin{equation}\label{eq:weighted-second-half}
 \begin{aligned}
 D_su&=\lim_{m\to\infty}Q_s^{1/2}\D(u_m(Y))
                         &&\text{in }L^2(\PP;H),\\
 \langle D_su,e_k\rangle_H&=q_k(s)^{1/2}d_k(u)
                         &&\text{in }L^2(\PP),\\
 \E_s(u,u)&=\EE\|D_su\|_H^2
             =\sum_{k\ge1}q_k(s)\|d_k(u)\|_2^2.
 \end{aligned}
\end{equation}
The last equality is Parseval's identity followed by Tonelli's
theorem. This proves \eqref{eq:weighted-domain}.
\end{proof}

\begin{proof}[Proof of Theorem~\ref{thm:Alimit}]
The form order follows directly from the weights in \eqref{eq:mosco-bridge}.
For $0<s<s'$ and $v\in\Ccore$, \eqref{eq:mosco-bridge} gives
\begin{equation}\label{eq:formorder-core}
 \E_{s'}(v,v)=\sum_kq_k(s')\|\partial_k(v(Y))\|_2^2
 \le\sum_kq_k(s)\|\partial_k(v(Y))\|_2^2
 =\E_s(v,v).
\end{equation}
An $\E_s+L^2$-Cauchy core sequence is consequently Cauchy in
the $\E_{s'}+L^2$ norm and has the same $L^2$ limit. Thus
\begin{equation}\label{eq:formorder}
 \F_s\subseteq\F_{s'},\qquad
 \E_{s'}(u,u)\le\E_s(u,u)\quad(u\in\F_s).
\end{equation}
Let $s_m\downarrow0$. For any $s'>0$ choose $m$ with $s_m<s'$.
Then \eqref{eq:formorder} implies
\begin{equation}\label{eq:limit-cofinality}
 \begin{aligned}
 \bigcap_m\F_{s_m}&=\bigcap_{s'>0}\F_{s'},\\
 \sup_m\E_{s_m}(u,u)&=\sup_{s'>0}\E_{s'}(u,u)
 \quad\left(u\in\bigcap_{s'>0}\F_{s'}\right).
 \end{aligned}
\end{equation}
Both the domain $\F_*$ and the energy $\E_*$ are therefore
independent of the chosen decreasing sequence.

The same core approximation, now using
$\E_s(v,v)\le\E_0^{\min}(v,v)$ on $\Ccore$, gives
\[
 \F_0^{\min}\subseteq\F_s,\qquad
 \E_s(u,u)\le\E_0^{\min}(u,u),\qquad u\in\F_0^{\min}.
\]
Thus the closed positive forms $\E_{s_m}$ increase and are bounded
above by the densely defined closed form $\E_0^{\min}$.
The monotone convergence theorem for forms
\cite[Theorem 3.1]{Simon} shows that $(\E_*,\F_*)$ is densely
defined and closed. Each normal contraction $\eta$ preserves all
$\F_{s_m}$ and satisfies
$\E_{s_m}(\eta(u),\eta(u))\le\E_{s_m}(u,u)$.
Taking the supremum proves the Markov property of $\E_*$.

For $v\in\Ccore$, monotone convergence in the Gaussian coordinates gives
\begin{equation}\label{eq:star-core-energy}
 \E_*(v,v)
 =\lim_m\sum_kq_k(s_m)\|\partial_k(v(Y))\|_2^2
 =\EE\|\D(v(Y))\|_H^2
 =\E_0^{\min}(v,v).
\end{equation}
An $\E_0^{\min}$-core approximation is therefore Cauchy in the
$\E_*+L^2(\rho)$ norm. Closedness yields
\begin{equation}\label{eq:star-min-inclusion}
 \F_0^{\min}\subseteq\F_*,\qquad
 \E_*(u,u)=\E_0^{\min}(u,u),\qquad u\in\F_0^{\min}.
\end{equation}
By Proposition~\ref{prop:forms}, $1\in\F_0^{\min}$ has zero
energy, so the limit form is conservative.

For the Mosco lower bound, let $s_l\downarrow0$ and
$u_l\rightharpoonup u$ in $L^2(\rho)$. Fix $s'>0$. Eventually
$s_l<s'$, and weak lower semicontinuity gives
\[
 \E_{s'}(u,u)\le\liminf_l\E_{s'}(u_l,u_l)
               \le\liminf_l\E_{s_l}(u_l,u_l).
\]
Taking the supremum over $s'>0$ proves
\begin{equation}\label{eq:star-mosco-liminf}
 \E_*(u,u)\le\liminf_l\E_{s_l}(u_l,u_l).
\end{equation}
For the upper bound, the constant sequence $u_l=u$ suffices:
$\E_{s_l}(u,u)\uparrow\E_*(u,u)$ when $u\in\F_*$, while the
bound is automatic when $\E_*(u,u)=+\infty$.
This proves Mosco convergence.

For $u\in\F_*$, Lemma~\ref{lem:weighted-gradient-domain}
provides the same weak derivatives $d_k(u)$ for each $s_l$.
Using \eqref{eq:weighted-domain} and monotone convergence,
\begin{equation}\label{eq:star-full-gradient-norm}
 \begin{aligned}
 \sum_{k\ge1}\|d_k(u)\|_2^2
 &=\lim_{l\to\infty}\sum_{k\ge1}q_k(s_l)
                       \|d_k(u)\|_2^2\\
 &=\sup_l\E_{s_l}(u,u)=\E_*(u,u)<\infty.
 \end{aligned}
\end{equation}
Set $F=u(Y)$ and $G=\sum_kd_k(u)e_k$.
By \eqref{eq:star-full-gradient-norm}, $G\in L^2(\PP;H)$.
Lemma~\ref{lem:weighted-gradient-domain} identifies
$\langle G,e_k\rangle_H$ as the weak coordinate derivative of
$F$. The finite-coordinate approximation used in the proof of
Lemma~\ref{lem:MalliavinY}, now with $p=2$, therefore gives
$F\in\mathbb D^{1,2}$ and $\D F=G$. Consequently
\begin{equation}\label{eq:star-sobolev-identification}
 \EE\|\D F\|_H^2
 =\sum_{k\ge1}\|d_k(u)\|_2^2=\E_*(u,u).
\end{equation}
Together with \eqref{eq:star-min-inclusion}, this proves
\eqref{eq:Agrad}. To prove strong locality, let $u,v\in\F_*$,
with $v=0$ $\rho$-almost everywhere outside a compact set
$K\subset E$ and $u$ constant $\rho$-almost everywhere on an
open neighborhood $U$ of $K$. Locality of the Gaussian Sobolev
gradient gives $\D(u(Y))=0$ on $\{Y\in U\}$ and
$\D(v(Y))=0$ on $\{Y\notin K\}$, almost surely. By polarization
of \eqref{eq:star-sobolev-identification},
\[
 \E_*(u,v)=\EE\langle\D(u(Y)),\D(v(Y))\rangle_H=0.
\]
Thus $\E_*$ is strongly local, completing Theorem~\ref{thm:Alimit}.
\end{proof}

\subsection{The endpoint integration-by-parts identity}
\begin{lemma}\label{lem:appendix-endpoint-ibp}
For $u\in C_c^\infty(\R^N)$ and
$g\in\F_*\cap L^\infty(\rho)$,
\begin{equation}\label{eq:appendix-endpoint-ibp}
 \E_*(u,g)=-\int_Eg\ell u\,d\rho.
\end{equation}
\end{lemma}
\begin{proof}
For every $s>0$, the domain inclusion $\F_*\subset\F_s$ and
Proposition~\ref{prop:IBP} give
\[
 \E_s(u,g)=-\int_Eg\ell_su\,d\rho.
\]
By the definition of the increasing limit and polarization,
$\E_s(u,g)\to\E_*(u,g)$. Proposition~\ref{prop:Acoeff} gives
$a_s\to a$ and $b_s\to b$ in $L^1(\rho)$, so
$\ell_su\to\ell u$ in $L^1(\rho)$. Since $g$ is bounded,
passing to the limit proves \eqref{eq:appendix-endpoint-ibp}.
\end{proof}

\section{A time-chain rule}\label{app:time-chain}
\begin{lemma}\label{lem:time-chain}
Let $(\E,\F)$ be a fixed response form with energy density
$\Gamma$, and let $\F'$ be the dual of $\F$ with pivot $L^2(\rho)$.
Suppose
\[
 f\in L^2(0,T;\F)\cap C([0,T];L^2(\rho)),\qquad
 \partial_tf\in L^2(0,T;\F'),\qquad c\le f\le C.
\]
Then $f^2\in L^2(0,T;\F)$ and
\begin{equation}\label{eq:time-chain-statement}
 \int_E\frac{f_b^3-f_a^3}{3}\,d\rho
 =\int_a^b\langle\partial_tf_t,f_t^2\rangle_{\F',\F}\,dt,
 \qquad 0\le a<b\le T.
\end{equation}
\end{lemma}
\begin{proof}
The form chain rule gives
$\E(f_t^2,f_t^2)\le4C^2\E(f_t,f_t)$, so the dual pairing is
integrable.
For $0<a<b<T$ and $0<\varepsilon<T-b$, set
$f_t^\varepsilon=\varepsilon^{-1}\int_t^{t+\varepsilon}f_r\,dr$.
Its derivative
$\partial_tf_t^\varepsilon=(f_{t+\varepsilon}-f_t)/\varepsilon$
belongs to $L^2(a,b;\F)$ and is the Steklov average of
$\partial_tf$ in $\F'$.
Thus $f^\varepsilon\in W^{1,2}(a,b;\F)$ and
$c\le f^\varepsilon\le C$.
Translation continuity gives $f^\varepsilon\to f$ in
$L^2(a,b;\F)$ and $\partial_tf^\varepsilon\to\partial_tf$
in $L^2(a,b;\F')$.

To pass the nonlinear test function to the limit, put
$w^\varepsilon=f^\varepsilon-f$.
Writing $(f^\varepsilon)^2-f^2=2fw^\varepsilon+(w^\varepsilon)^2$,
the energy product rule yields
\begin{equation}\label{eq:chain-energy}
 \begin{aligned}
 \Gamma((f_t^\varepsilon)^2-f_t^2,
                         (f_t^\varepsilon)^2-f_t^2)
 &\le8C^2\Gamma(w_t^\varepsilon,w_t^\varepsilon)\\
 &\quad+8|w_t^\varepsilon|^2\Gamma(f_t,f_t).
 \end{aligned}
\end{equation}
After integration in $t,\rho$, the first term tends to zero
by strong form convergence.
For the second, $w^\varepsilon\to0$ in $dt\,d\rho$-measure
and $|w^\varepsilon|\le2C$, while the fixed weight
$\Gamma(f,f)$ is integrable.
Together with $|(f^\varepsilon)^2-f^2|\le2C|w^\varepsilon|$,
this proves $(f^\varepsilon)^2\to f^2$ in $L^2(a,b;\F)$.

Apply the $L^2(\rho)$ chain rule to a smooth extension of
$r^3/3$ from $[c,C]$ with bounded second derivative.
Since $f^\varepsilon$ stays in $[c,C]$, it gives
\begin{equation}\label{eq:stechain}
 \frac{d}{dt}\int_E\frac{(f_t^\varepsilon)^3}{3}\,d\rho
 =\langle\partial_tf_t^\varepsilon,
                     (f_t^\varepsilon)^2\rangle_{\F',\F}.
\end{equation}
Integrate over $(a,b)$. The strong convergences pass the dual
pairing to the limit by Cauchy--Schwarz.
The endpoint terms converge by $L^2(\rho)$
continuity and $|r^3-s^3|\le3C^2|r-s|$ on $[c,C]$.
This proves \eqref{eq:time-chain-statement}, including $a=0$ and
$b=T$ by continuity. Integrability of the pairing also gives
absolute continuity of the cubic energy.
\end{proof}

\end{document}